\documentclass[a4paper,12pt]{article}
\usepackage{amssymb,fullpage,amsmath,amsthm}

\newtheorem{thm}{Theorem}

\newtheorem{lem}[thm]{Lemma}
\newtheorem{propn}[thm]{Proposition}

\newtheorem{assu}{Assumption}

\newenvironment{assu1hat}[1][Assumption \^{1}]{ \begin{trivlist}\item[\hskip\labelsep{\bfseries #1}]}{\end{trivlist}}
\newenvironment{assu2hat}[1][Assumption \^{2}]{ \begin{trivlist}\item[\hskip\labelsep{\bfseries #1}]}{\end{trivlist}}

\newenvironment{assu1r}[1][Assumption {1}R]{ \begin{trivlist}\item[\hskip\labelsep{\bfseries #1}]}{\end{trivlist}}
\newenvironment{assu2r}[1][Assumption {2}R]{ \begin{trivlist}\item[\hskip\labelsep{\bfseries #1}]}{\end{trivlist}}
\newenvironment{assu3r}[1][Assumption {3}R]{ \begin{trivlist}\item[\hskip\labelsep{\bfseries #1}]}{\end{trivlist}}

\begin{document}

\title{Local limit theorems for sequences of simple random walks on graphs}

\author{D.A. Croydon\footnote{Dept of Statistics,
University of Warwick, Coventry CV4 7AL, UK;
{d.a.croydon@warwick.ac.uk}.} and B. M. Hambly\footnote{Mathematical Institute, University of Oxford, 24-29 St Giles', Oxford, OX1 3LB, UK; hambly@maths.ox.ac.uk}}
\date{1 October 2008}

\maketitle

\begin{abstract}
In this article, local limit theorems for sequences of simple random walks on graphs are established. The results formulated are motivated by a variety of random graph models, and explanations are provided as to how they apply to supercritical percolation clusters, graph trees converging to the continuum random tree and the homogenisation problem for nested fractals. A subsequential local limit theorem for the simple random walks on generalised Sierpinski carpet graphs is also presented.
\end{abstract}

\section{Introduction}
The classical local limit theorem (see \cite{Feller}, Section XV.5, for example) describes how the transition probabilities of the discrete time simple random walk on $\mathbb{Z}$ can be rescaled to yield the Gaussian transition densities of Brownian motion on $\mathbb{R}$. Analogous statements have been proved in many other settings, including the recent result of \cite{BarHamllt}, which demonstrates that the transition probabilities of the discrete time simple random walk on the random environment generated by supercritical bond-percolation on $\mathbb{Z}^d$ can be rescaled to a Gaussian limit (for almost-every environment). In this article, by generalising the argument of \cite{BarHamllt}, we deduce that the corresponding limit result holds for any sequence of simple random walks on graphs whose laws can be rescaled appropriately and which satisfies a tightness assumption on its transition densities.

Let us start by describing the framework of this article. We assume that there exists an underlying metric space $(E,d_E)$ and suppose $F$ is a subset of $E$ such that $F\cap \overline{B}_E(x,r)$ is compact for every $x\in E$ and $r>0$ (where $\overline{B}_E(x,r)$ is the closed ball in $(E,d_E)$ with centre $x$ and radius $r$). We also presume that the following are defined: $\rho$, a distinguished element of $F$; $\nu$, a Radon measure with full support on the metric space $(F,d_F)$, where $d_F:=d_E|_{F\times F}$; and $(q_t(x))_{x\in F, t>0}$, a family of densities so that, for each $t>0$, $q_t$ is a Borel measurable non-negative function on $F$ which integrates to 1 with respect to $\nu$. Moreover, we suppose that $(q_t(x))$ is a jointly continuous function of $(t,x)$. Typically in examples we have a conservative $\nu$-symmetric Markov diffusion $X$ with transition density $(p_t(x,y))_{x,y\in F, t>0}$, and in this case we take $q_t(x)=p_t(\rho,x)$. The entities introduced above will represent the limits of sequences of corresponding objects defined from sequences of graphs. Note that our assumptions imply that $F$ is a closed, and therefore measurable, subset of $E$, therefore the measure $\nu$ can be extended to a Borel measure on $(E,d_E)$.

We continue by introducing some general notation for random walks on graphs. First, fix $G=(V(G),E(G))$ to be a locally finite connected graph with at least two vertices, where $V(G)$ denotes the vertex set and $E(G)$ the edge set of $G$. For $x,y\in V(G)$, we write the number of edges in the shortest path from $x$ to $y$ in $G$ as $d_G(x,y)$, so that $d_G$ is a metric on $V(G)$. Define a symmetric weight function $\mu^G:V(G)^2\rightarrow \mathbb{R}_+$ that satisfies $\mu^G_{xy}>0$ if and only if $\{x,y\}\in E(G)$. The discrete time simple random walk on the weighted graph $G$ is then the Markov chain $((X^G_m)_{m\geq 0}, \mathbf{P}^G_x,x\in V(G))$ with transition probabilities $(P_G(x,y))_{x,y\in V(G)}$ defined by $P_G(x,y):={\mu^G_{xy}}/{\mu^G_{x}}$, where $\mu^G_x:=\sum_{y\in V(G)}\mu^G_{xy}$. If we define a measure $\nu^G$ on $V(G)$ by setting, for $A\subseteq V(G)$, $\nu^G(A):=\sum_{x\in A}\mu^G_x$,
then $\nu^G$ is invariant for $X^G$, and the transition density of $X^G$, with respect to $\nu^G$, is given by $(p^G_m(x,y))_{x,y\in V(G),\:m\geq 0}$, where
\[p_m^G(x,y):=\frac{\mathbf{P}^{G}_x(X_m=y)}{\nu^G(\{y\})}.\]
Due to parity concerns for bipartite graphs, to obtain a smooth limiting result, rather than the transition density itself, we will consider $(q^G_m(x,y))_{x,y\in V(G),\:m\geq 0}$ defined by
\begin{equation}\label{qdef}
q^G_m(x,y):=\frac{p^G_m(x,y)+p^G_{m+1}(x,y)}{2},
\end{equation}
and also define $q_m^G(x):=q_m^G(\rho(G),x)$, where $\rho(G)$ is a distinguished element of $V(G)$.

For our main local limit theorem, we suppose that a sequence of graphs $(G^n)_{n\geq 1}$ have been embedded into $E$ so that the various sequences of objects described in the previous paragraph approximate $d_E$, $F$, $\nu$ and the laws associated with $q_t$, $t>0$, in the way we now describe precisely. For brevity, throughout this article we write $\nu^{G^n}, X^{G^n}, q^{G^n},\dots$ as $\nu^n, X^n, q^n,\dots$ respectively.

\begin{assu} \label{first} Let $(G^n)_{n\geq 1}$ be a sequence of locally finite connected graphs that satisfy $\# V(G^n)\geq 2$, $V(G^n)\subseteq E$ and $\rho(G^n)=\rho$ for every $n\geq 1$. Fix three non-negative divergent sequences $(\alpha(n))_{n\geq 1}$, $(\beta(n))_{n\geq 1}$, and $(\gamma(n))_{n\geq 1}$, and suppose that:
 \renewcommand{\labelenumi}{(\alph{enumi})}
\begin{enumerate}
  \item there exists a constant $c_1>0$ such that, for $n\geq 1$,
  \begin{equation}\label{as1a}
  d_{G^n}(x,y)\geq c_1\alpha(n)d_E(x,y),\hspace{20pt}\forall x,y\in V(G^n).
  \end{equation}
  Furthermore, there exists a non-negative sequence $(\tilde{\alpha}(n))_{n\geq 1}$, such that $\tilde{\alpha}(n)=o(\alpha(n))$ as $n\rightarrow \infty$ and also, for each $r>0$, there exists a finite constant $c_2$ and an integer $n_0$ such that
  \[d_{G^n}(x,y)\leq c_2\alpha(n)d_E(x,y)+\tilde{\alpha}(n),\hspace{20pt}\forall x,y\in V(G^n)\cap B_E(\rho,r),\]
  for $n\geq n_0$, where $B_E(\rho,r)$ is the open ball in $(E,d_E)$ with centre $\rho$ and radius $r$.
  \item for each $r>0$,
  \[\lim_{n\rightarrow\infty}\sup_{x\in B_F(\rho,r)}d_E(x,V(G^n))= 0,\]
  where $B_F(\rho,r)$ is the open ball in $(F,d_F)$ with centre $\rho$ and radius $r$.
  \item for every $x\in F$ and $r>0$,
  \[\lim_{n\rightarrow\infty}\beta(n)^{-1}\nu^{n}(B_E(x,r))= \nu(B_E(x,r)).\]
  \item for any compact interval $I\subset (0,\infty)$, $x\in F$ and $r>0$,
\[\lim_{n\rightarrow\infty}\mathbf{P}_{\rho}^{G^n}\left(X^n_{\lfloor \gamma(n)t\rfloor}\in B_E(x,r)\right)=\int_{B_F(x,r)}q_t(y)\nu(dy)\]
uniformly for $t\in I$.
\end{enumerate}
\end{assu}

In addition to these approximation conditions, we will apply the following tightness condition for the transition densities of the simple random walks on the graphs $(G^n)_{n\geq1}$. In the case when $V(G^n)\subseteq F$ for every $n$, we will show that (given Assumption 1) this condition is actually necessary for a local limit theorem of the type we prove to hold.

\begin{assu}\label{second} In the setting of Assumption \ref{first}, suppose that, for any compact interval $I\subset (0,\infty)$ and $r>0$,
\[\lim_{\delta\rightarrow 0}\limsup_{n\rightarrow\infty}\sup_{\substack{x,y\in B_{G^n}(\rho,\alpha(n)r):\\d_{G^n}(x,y)\leq \alpha(n)\delta}}\sup_{t \in I}\beta(n)\left|q_{\lfloor\gamma(n)t\rfloor}^{n}(x)-q_{\lfloor\gamma(n)t\rfloor}^{n}(y)\right|=0,\]
where $B_{G^n}(\rho,r)$ is the open ball in $(V(G^n),d_{G^n})$ with centre $\rho$ and radius $r$.
\end{assu}

Finally, before we state our first main result, observe that if Assumption \ref{first} holds, then for $r>0$, $n\geq 1$, we can bound the graph distance $d_{G^n}(x,y)$ above by a constant (depending on $n$) uniformly over $x,y\in V(G^n)\cap B_E(\rho,r)$. Hence, because $G^n$ is by definition a locally finite graph, there can only be a finite number of points in the set $V(G^n)\cap B_E(\rho,r)$, and consequently the same is true for any set of the form $V(G^n)\cap B_E(x,r)$. In particular, this implies that for every $x\in E$, we can choose (not necessarily uniquely) a point $g_n(x)\in V(G^n)$ that minimises the distance $d_E(x,y)$ over $y\in V(G^n)$.

{\thm \label{bdd} Fix a compact interval $I\subset (0,\infty)$ and $r>0$. Suppose Assumptions \ref{first} and \ref{second} hold, then
\begin{equation}\label{conc}
\lim_{n\rightarrow\infty}\sup_{x\in B_F(\rho,r)}\sup_{t\in I}\left|\beta(n)q^n_{\lfloor\gamma(n)t\rfloor}(g_n(x))-q_t(x)\right|=0.
\end{equation}
Conversely, when $V(G^n)\subseteq F$ for every $n$, if Assumption 1 is satisfied and (\ref{conc}) holds for every compact interval $I\subset (0,\infty)$ and $r>0$, then Assumption 2 holds.}
\bigskip

To allow us to extend (\ref{conc}) to hold uniformly over unbounded time intervals and non-compact spaces $F$, we need to impose some extra conditions which guarantee the decay in time and space of the discrete and continuous transition densities.

\begin{assu}\label{third} The metric space $(F,d_F)$ has the midpoint property, i.e. for every $x,y\in F$, there exists $z\in F$ such that $d_F(x,z)=\frac{1}{2}d_F(x,y)=d_F(z,y)$. Furthermore, the following conditions are fulfilled.
 \renewcommand{\labelenumi}{(\alph{enumi})}
\begin{enumerate}
\item The transition density of $X$ satisfies
\[\lim_{t\rightarrow\infty}\sup_{x\in F} q_t(x)=0,\hspace{20pt}\lim_{r\rightarrow\infty}\sup_{x\in F\backslash B_F(\rho,r)}\sup_{t\in I}q_t(x)=0,\]
for any compact interval $I\subset (0,\infty)$.
\item In the setting of Assumption \ref{first},
\[\lim_{t\rightarrow\infty}\limsup_{n\rightarrow\infty}\sup_{x\in V(G^n)} \beta(n)q_{\lfloor \gamma(n)t\rfloor}^n(x)=0,\]
and, for any compact interval $I\subset (0,\infty)$,
\[\lim_{r\rightarrow\infty}\limsup_{n\rightarrow\infty}\sup_{x\in V(G^n)\backslash B_{G^n}(\rho,\alpha(n)r)} \sup_{t\in I}\beta(n)q_{\lfloor \gamma(n)t\rfloor}^n(x)=0.\]
\end{enumerate}
\end{assu}

If this extra assumption is satisfied, then we are able to prove the following.

{\thm \label{unbdd} Fix $T_1>0$. Suppose Assumptions \ref{first}, \ref{second} and \ref{third} hold, then
\[\lim_{n\rightarrow\infty}\sup_{x\in F}\sup_{t\geq T_1}\left|\beta(n)q^n_{\lfloor\gamma(n)t\rfloor}(g_n(x))-q_t(x)\right|=0.\]}

The main motivation for proving results such as Theorems \ref{bdd} and \ref{unbdd} is to provide conditions under which a weak convergence result, such as that appearing in Assumption \ref{first}(d), implies a local limit theorem. Obviously, the usefulness of such results depends on the applicability of the assumptions that have been made, and so, after completing the proofs of Theorems \ref{bdd} and \ref{unbdd} in Section \ref{proof}, we provide two alternative sufficient conditions for Assumption \ref{second}. The first of these, see Assumption \ref{fourth} in Section \ref{PHI}, involves the parabolic Harnack inequality, which is also known to imply various other analytic conditions for random walks on graphs (see \cite{Kumagai} for a summary of results in this area). The second, see Assumption \ref{fifth} in Section \ref{ressec}, relies on being able to bound the resistance metric on graphs in the sequence $(G^n)_{n\geq 1}$ using the shortest path metric, and, as we shall see in Section \ref{examples}, is applicable to graph trees and nested fractal graphs.

In Section \ref{2param}, a short investigation into the asymptotics of $(q^n_m(x,y))_{x,y\in V(G),m\geq0}$, when considered as a function of two spatial coordinates, is presented. In particular, we give sufficient conditions for the uniform convergence of
\[(\beta(n)q^n_{\lfloor \gamma(n) t\rfloor}(g_n(x),g_n(y)))_{x,y\in F,t>0}\]
to the transition density of a Markov process $(p_t(x,y))_{x,y\in F, t>0}$ (at least in bounded space-time regions). This is followed, in Section \ref{randsec}, by a demonstration of how the analytic results that we prove can be adapted to the case when the weights of the graphs $(G^n)_{n\geq 1}$ are random and we only have probabilistic versions of Assumptions 1, 2 and 3 instead of almost-sure versions. We conclude our article with a collection of examples to which our results apply, including random graphs on the integer lattice, graph trees converging to the continuum random tree, nested fractal graphs and generalised Sierpinski carpet graphs.

Finally, define the continuous time simple random walk on a graph $G$ to be the continuous time Markov process $((Y_t^G)_{t\geq 0},\tilde{\mathbf{P}}_x^G,x\in V(G))$ with generator $\mathcal{L}_G$, as defined below at (\ref{generator}). The transition density of $Y^G$, with respect to its invariant measure $\nu^G$, is given by
\[\tilde{p}_t^G(x,y)=\frac{\tilde{\mathbf{P}}_x^G(Y^G_t=y)}{\nu^G(\{y\})},\]
and we write $\tilde{q}^G_t(x):=\tilde{p}_t^G(\rho,x)$. Under the continuous time analogues of Assumptions 1, 2 and 3, it is possible to obtain continuous time versions of Theorems \ref{bdd} and \ref{unbdd} that apply to the continuous time simple random walks on the graphs $(G^n)_{n\geq 1}$. However, since they can be proved using identical arguments, we omit them. See Examples \ref{nf}, \ref{vicsek} and \ref{carpetsec} for results which illustrate the more general continuous time local limit theorems.

\section{Proof of local limit theorems}\label{proof}

The aim of this section is to prove Theorems \ref{bdd} and \ref{unbdd}. We start by generalising slightly Assumption \ref{first}(d).

{\lem\label{gen} Suppose Assumption \ref{first}(d) holds, then, for any compact interval $I\subset (0,\infty)$, $x\in F$ and $r>0$,
\[\lim_{n\rightarrow\infty}\mathbf{P}_{\rho}^{G^n}\left(X^n_{\lfloor \gamma(n)t\rfloor+i}\in B_E(x,r)\right)=\int_{B_F(x,r)}q_t(y)\nu(dy)\]
uniformly for $t\in I$, $i\in \{0,1\}$.}
\begin{proof} The proof is elementary, and requires the application of only Assumption \ref{first}(d), the joint continuity of $(q_t(x))$ in $(t,x)$ and the fact that $F\cap\overline{B}_E(x,r)$ is compact.
\qed\end{proof}

We now prove a point-wise version of a local limit theorem.

{\propn\label{point} Fix a compact interval $I\subset (0,\infty)$ and suppose Assumptions \ref{first} and \ref{second} hold, then, for every $x\in F$,
\[\lim_{n\rightarrow\infty}\sup_{t\in I}\left|\beta(n)q^n_{\lfloor\gamma(n)t\rfloor}(g_n(x))-q_t(x)\right|=0.\]}
\begin{proof} Fix $x\in F$, $\varepsilon>0$ and set $r=d_E(\rho,x)$. Let $c$ be a finite constant and $n_0$ an integer such that $d_{G^n}(y,z)\leq c\alpha(n)d_E(y,z)+\tilde{\alpha}(n)$ for every $y,z\in V(G^n)\cap B_E(\rho,r+1)$ and $n\geq n_0$; the existence of such constants is guaranteed by Assumption \ref{first}(a). Furthermore, by the tightness condition of Assumption \ref{second}, and the supposition that $(q_t(y))$ is jointly continuous in $(t,y)$, we can choose $r_0\in(0,1)$ small enough and an integer $n_1\geq n_0$ so that
\begin{equation}\label{dischkcont}
\sup_{\substack{y,z\in B_{G^n}(\rho,c\alpha(n)(r+2)):\\d_{G^n}(y,z)\leq 3c\alpha(n)r_0}}\sup_{t\in I}\beta(n)\left|q_{\lfloor\gamma(n)t\rfloor}^{n}(y)-q_{\lfloor\gamma(n)t\rfloor}^{n}(z)\right|\leq \varepsilon,
\end{equation}
for $n\geq n_1$, and also
\begin{equation}\label{hkcont}
\sup_{\substack{y,z\in B_F(\rho,r+1):\\d_F(y,z)\leq r_0}}\sup_{t\in I}|q_t(y)-q_t(z)|\leq \varepsilon.
\end{equation}

For this choice of $r_0$, we consider the quantity
\begin{eqnarray*}
\lefteqn{J(t,n):=\frac{1}{2}\left\{\mathbf{P}_{\rho}^{G^n}\left(X^n_{\lfloor \gamma(n)t\rfloor}\in B_E(x,r_0)\right)+\mathbf{P}_{\rho}^{G^n}\left(X^n_{\lfloor \gamma(n)t\rfloor+1}\in B_E(x,r_0)\right)\right\}}\\
&&\hspace{240pt}-\int_{B_F(x,r_0)}q_t(y)\nu(dy),
\end{eqnarray*}
which can be written as $J(t,n)=J_1(t,n)+J_2(t,n)+J_3(t,n)+J_4(t,n)$, where
\begin{eqnarray*}
J_1(t,n)&:=&\sum_{y\in V(G^n)\cap B_E(x,r_0)} \left(q^n_{\lfloor \gamma(n)t\rfloor} (y)-q^n_{\lfloor \gamma(n)t\rfloor}(g_n(x))\right)\nu^n(\{y\}),\\
J_2(t,n)&:=&\beta(n)^{-1}\nu^n(B_E(x,r_0))\left(\beta(n)q^n_{\lfloor \gamma(n)t\rfloor}(g_n(x))- q_t(x)\right),\\
J_3(t,n)&:=&q_t(x)\left(\beta(n)^{-1}\nu^n(B_E(x,r_0))-\nu(B_E(x,r_0))\right),\\
J_4(t,n)&:=&\int_{B_F(x,r_0)}\left(q_t(x)-q_t(y)\right)\nu(dy).
\end{eqnarray*}

Now, by (\ref{hkcont}), we immediately have $\sup_{t\in I}|J_4(t,n)|\leq \varepsilon \nu(B_E(x,r_0))$. Furthermore, by applying Assumption \ref{first}(c) and Lemma \ref{gen}, it is possible to deduce that there exists an integer $n_2\geq n_1$ such that $\sup_{t\in I}|J_3(t,n)|,\sup_{t\in I}|J(t,n)|\leq \varepsilon \nu(B_E(x,r_0))$ for $n\geq n_2$. To bound $J_1(t,n)$ in a similar fashion, first note that Assumption \ref{first}(b) implies $d_E(g_n(x),x)\rightarrow 0$ as $n\rightarrow \infty$. In particular, it follows that there exists an integer $n_3\geq n_2$ such that $g_n(x)\in B_E(x,r_0)$ for $n\geq n_3$. Thus, for $n\geq n_3$,
\begin{equation}\label{j1bound}
\sup_{t\in I}|J_1(t,n)|\leq \sup_{y,z\in V(G^n)\cap B_E(x,r_0)}\sup_{t\in I}\left|
q^n_{\lfloor \gamma(n)t\rfloor} (y)-q^n_{\lfloor \gamma(n)t\rfloor}(z)\right|\nu^n(B_E(x,r_0)).
\end{equation}
Recall from Assumption \ref{first}(a) that $\tilde{\alpha}(n)=o(\alpha(n))$, and therefore we can choose an integer $n_4\geq n_3$ such that $\tilde{\alpha}(n)\leq c\alpha(n)r_0$. Consequently, for $y,z\in V(G^n)\cap B_E(x,r_0)$, we have that $y,z\in B_{G^n}(\rho,c\alpha(n) (r+2))$ and also $d_{G^n}(y,z)\leq 3c\alpha(n)r_0$ whenever $n\geq n_4$. Hence, by (\ref{dischkcont}),
\[\sup_{t\in I}|J_1(t,n)|\leq \varepsilon \beta(n)^{-1}\nu^n(B_E(x,r_0))\]
for $n\geq n_4$. To bound the right-hand side of this expression, note that Assumption \ref{first}(c) allows us to choose $n_5\geq n_4$ such that $\left|\beta(n)^{-1}\nu^n(B_E(x,r_0))-\nu(B_E(x,r_0))\right|\leq \nu(B_E(x,r_0))/2$ for $n\geq n_5$. Thus $\sup_{t\in I}|J_1(t,n)|\leq 2\varepsilon\nu (B_E(x,r_0))$ for $n\geq n_5$. Piecing all these bounds together yields, for $n\geq n_5$,
$\sup_{t\in I} |J_2(t,n)|\leq 5\varepsilon \nu(B_E(x,r_0))$. Finally, note that the left-hand side of this expression is bounded below by
\[\sup_{t\in I}\left|\beta(n)q^n_{\lfloor \gamma(n)t\rfloor}(g_n(x)) - q_t(x)\right| \nu(B_E(x,r_0))/2\]
whenever $n\geq n_5$. The result follows.
\qed\end{proof}

This result is readily extended to hold uniformly over bounded balls in $F$, thereby establishing Theorem \ref{bdd}.

\begin{proof}[Proof of Theorem \ref{bdd}] Suppose Assumptions 1 and 2 hold. Fix $r, \varepsilon>0$ and choose $c$, $r_0$ and $n_1$ as in the proof of the previous proposition. By assumption, $\overline{B}_F(\rho,r)$ is compact, hence there exists a finite collection $\mathcal{X}\subseteq \overline{B}_F(\rho,r)$ such that $(B_F(x,r_0))_{x\in\mathcal{X}}$ is an open cover for $B_F(\rho,r)$. Since $\mathcal{X}$ is finite, applying Proposition \ref{point} allows it to be deduced that there exists an integer $n_2\geq n_1$ such that
\begin{equation}\label{finite}
\sup_{x\in\mathcal{X}}\sup_{t\in I}\left|\beta(n)q^n_{\lfloor\gamma(n)t\rfloor}(g_n(x))-q_t(x)\right|\leq \varepsilon,
\end{equation}
for $n\geq n_2$. Now, suppose $x\in B_F(\rho,r)$, then $x\in B_F(y(x),r_0)$ for some $y(x)\in \mathcal{X}$, and we can write
\begin{eqnarray*}
\sup_{t\in I} \left|\beta(n)q^n_{\lfloor\gamma(n)t\rfloor}(g_n(x))-q_t(x)\right|&\leq&\sup_{t\in I} \beta(n)\left|q^n_{\lfloor\gamma(n)t\rfloor}(g_n(x))-q^n_{\lfloor\gamma(n)t\rfloor}(g_n(y(x)))\right|\\
&&+\sup_{t\in I} \left|\beta(n)q^n_{\lfloor\gamma(n)t\rfloor}(g_n(y(x)))-q_t(y(x))\right|\\
&&+\sup_{t\in I} \left|q_t(y(x))-q_t(x)\right|.
\end{eqnarray*}
Since $x,y(x)\in B_E(\rho, r+1)$ and $d_E(x,y(x))\leq r_0$, the inequality at (\ref{hkcont}) implies that the final term here is bounded above by $\varepsilon$ uniformly over $x\in B_F(\rho,r)$. It follows from (\ref{finite}) that the second term is also bounded above by $\varepsilon$ uniformly over $x\in B_F(\rho,r)$  for $n\geq n_2$. To deal with the first term, we start by choosing an integer $n_3\geq n_2$ such that $\tilde{\alpha}(n)\leq c\alpha(n)r_0$ and $d_E(x,g_n(x))<r_0$ for every $x\in\overline{B}_F(\rho,r)$ and $n\geq n_3$ (this is possible by Assumptions \ref{first}(a) and \ref{first}(b)). For $n\geq n_3$, we therefore have $g_n(x),g_n(y(x))\in V(G^n)\cap B_E(\rho,r+1)$, and consequently, as in the proof of the previous proposition, $g_n(x),g_n(y(x))\in B_{G^n}(\rho,c\alpha(n) (r+2))$ and also $d_{G^n}(g_n(x),g_n(y(x)))\leq 3c\alpha(n)r_0$ whenever $n\geq n_3$. Thus
we can apply the inequality at (\ref{dischkcont}) to deduce that the first term in the above upper bound is also bounded by $\varepsilon$ uniformly over $x\in B_F(\rho,r)$, and the proof of (\ref{conc}) is complete.

Let us now assume that $V(G^n)\subseteq F$, Assumption 1 is satisfied and (\ref{conc}) holds for every compact interval $I\subset (0,\infty)$ and $r>0$. Setting $c_1$ to be the constant of (\ref{as1a}), it is clear that
\begin{eqnarray*}
\lefteqn{\sup_{\substack{x,y\in B_{G^n}(\rho,\alpha(n)r):\\d_{G^n}(x,y)\leq \alpha(n)\delta}}\sup_{t \in I}\beta(n)\left|q_{\lfloor\gamma(n)t\rfloor}^{n}(x)-q_{\lfloor\gamma(n)t\rfloor}^{n}(y)\right|}\\
&\leq&2\sup_{x\in B_F(\rho,c_1^{-1}r)}\sup_{t\in I}\left|\beta(n)q^n_{\lfloor\gamma(n)t\rfloor}(g_n(x))-q_t(x)\right|\\
&&+\sup_{\substack{x,y\in B_{F}(\rho,c_1^{-1}r):\\d_{F}(x,y)\leq c_1^{-1}\delta}}\sup_{t \in I}\left|q_{t}(x)-q_{t}(y)\right|
\end{eqnarray*}
Assumption 2 is readily deduced from this bound by applying the limit at (\ref{conc}) and the joint continuity of $(q_t(x))_{x\in F,t>0}$.
\qed\end{proof}

To complete this section, we demonstrate how Assumption \ref{third} allows this result to be extended to unbounded regions of time and space. However, we first prove a simple lemma relating the $d_F$-distance of a point $x$ from $\rho$ in $F$ to the $d_{G^n}$-distance of the point $g_n(x)$ from $\rho$ in $G^n$ whenever the midpoint property is satisfied by $(F,d_F)$.

{\lem \label{lem} Suppose Assumption \ref{first} holds, let $c=c_1/2$, where $c_1$ is the constant of the bound at (\ref{as1a}), and assume that $(F,d_F)$ satisfies the midpoint property. If $r>0$, then there exists an integer $n_0$ such that
\[\inf_{x\in F\backslash B_F(\rho,r)}d_{G^n}(\rho,g_n(x))\geq c\alpha(n) r,\]
for every $n\geq n_0$.}
\begin{proof} Fix $r>0$. We first observe that the midpoint property of $(F,d_F)$ and the assumption that sets of the form $\overline{B}_F(x,r)$, $x\in F$, $r>0$, are compact imply that for each $x,y\in F$ there exists a (not necessarily unique) geodesic path $(\gamma(t))_{0\leq t\leq 1}$ in $F$ such that $\gamma(0)=x$, $\gamma(1)=y$ and $d_F(\gamma(s),\gamma(t))=|t-s|d_F(x,y)$ for $0\leq s\leq t\leq 1$. Hence, for every $x\in F\backslash B_F(\rho,r)$, there exists a $y(x)$ such that $d_F(\rho,y(x))=3r/4$ and $d_F(y(x),x)=d_F(\rho,x)-3r/4$. Consequently, we have that
\[d_E(\rho,y(x))+d_E(y(x),x) = d_E(\rho,x)\leq d_E(\rho,g_n(x))+d_E(g_n(y(x)),y(x))+d_E(y(x),x)\]
where we have applied the definition of $g_n(x)$ as the closest point in $V(G^n)$ to $x$. Canceling $d_E(y(x),x)$ from each side yields $3r/4\leq d_E(\rho,g_n(x))+d_E(g_n(y(x)),y(x))$.

By Assumption \ref{first}(b), we can choose an integer $n_0$ such that
\[\sup_{y\in B_F(\rho,r)}d_E(y,g_n(y))\leq r/4,\]
for $n\geq n_0$. Thus we can conclude from the previous paragraph that $d_E(\rho, g_n(x))\geq r/2$, for every $x\in F\backslash B_F(\rho,r)$ and $n\geq n_0$. The result is an easy consequence of this bound and Assumption \ref{first}(a).
\qed\end{proof}

\begin{proof}[Proof of Theorem \ref{unbdd}] Fix $T_1$ and $\varepsilon>0$. By Assumption \ref{third}, we can choose a finite time $T_2\geq T_1$ and an integer $n_0$ such that $\sup_{x\in V(G^n)} \beta(n)q_{\lfloor \gamma(n)t\rfloor}^n(\rho,x)\leq \varepsilon$, for $t\geq T_2$ and $n\geq n_0$, and also $\sup_{x\in F} q_t(x)\leq \varepsilon$ for $t\geq T_2$. Clearly, for this choice, we have
\begin{equation}\label{AAA}
\sup_{x\in F}\sup_{t\geq T_2}\left|\beta(n)q^n_{\lfloor\gamma(n)t\rfloor}(g_n(x))-q_t(x)\right|\leq 2\varepsilon
\end{equation}
for $n\geq n_0$.

Taking $I=[T_1,T_2]$, applying Assumption \ref{third} allows it to be deduced that there exists a finite radius $r_0$ and integer $n_1\geq n_0$ such that
\begin{equation}\label{near}
\sup_{x\in V(G^n)\backslash B_{G^n}(\rho,\alpha(n)r_0)} \sup_{t\in I}\beta(n)q_{\lfloor \gamma(n)t\rfloor}^n(x)\leq \varepsilon,
\end{equation}
for $n\geq n_1$, and also $\sup_{x\in F\backslash B_F(\rho,r_0)}\sup_{t\in I}q_t(x)\leq \varepsilon$. Now, let $c$ be the constant of Lemma \ref{lem} and define $r_1:=r_0(1+c^{-1})$. By Lemma \ref{lem}, there exists an integer $n_2\geq n_1$ such that for every $x\in F\backslash B_F(\rho,r_1)$ we have $g_n(x)\in V(G^n)\backslash B_{G^n}(\rho,\alpha(n)r_0)$ for $n\geq n_2$, and so we can apply the inequality at (\ref{near}) to deduce that
\[\sup_{x\in F\backslash B_{F}(\rho,r_1)} \sup_{t\in I}\beta(n)q_{\lfloor \gamma(n)t\rfloor}^n(g_n(x))\leq \varepsilon,\]
for $n\geq n_2$. Thus, because it also holds that $\sup_{x\in F\backslash B_F(\rho,r_1)}\sup_{t\in I}q_t(x)\leq \varepsilon$, it follows that
\[\sup_{x\in F\backslash B_F(\rho,r_1)}\sup_{t\in I}\left|\beta(n)q^n_{\lfloor\gamma(n)t\rfloor}(g_n(x))-q_t(x)\right|\leq 2\varepsilon,\]
for $n\geq n_2$. To complete the proof, it suffices to combine this conclusion with (\ref{AAA}) and the convergence result of Theorem \ref{bdd}.
\qed\end{proof}

\section{Parabolic Harnack inequality and tightness}\label{PHI}

For a locally finite connected graph $G$ define, for $x\in V(G)$, $R,T\geq 0$,
\[Q_G(x,R,T):=[0,T]\times B_G(x,R),\]
and also
\[Q^-_G(x,R,T):=[\tfrac{1}{4}T,\tfrac{1}{2}T]\times B_G(x,\tfrac{1}{2}R),\hspace{10pt}Q^+_G(x,R,T):=[\tfrac{3}{4}T,T)\times B_G(x,\tfrac{1}{2}R).\]
We describe a function $u(n,x)$ as caloric on $Q_G(x,R,T)$ if $u$ is defined on the set $\overline{Q}_G(x,R,T):=([0,T]\cap \mathbb{Z})\times {B}_G(x,R+1)$, and
\[u(n+1,x)-u(n,x)=\mathcal{L}_Gu(n,x),\]
for every $0\leq n\leq T-1$ and $x\in B_G(x,R)$, where $\mathcal{L}_G$ is the generator of the random walk $X^G$, which can be defined as the operator satisfying
\begin{equation}\label{generator}
\mathcal{L}_Gf(x)=\sum_{y\in V(G)}P_G(x,y)(f(y)-f(x)),
\end{equation}
for functions $f\in \mathbb{R}^{V(G)}$. The parabolic Harnack inequality with constant $C_H$ is then said to hold for $Q_G(x,R,T)$ if whenever $u$ is non-negative and caloric on $Q_G(x,R,T)$, we have
\[\sup_{(n,x)\in Q^-_G(x,R,T)} \hat{u}(n,x)\leq C_H\inf_{(n,x)\in Q^+_G(x,R,T)} \hat{u}(n,x),\]
where $\hat{u}(n,x):=u(n+1,x)+u(n,x)$.

We show in this section how if we assume that the sequence of graphs $(G^n)_{n\geq 1}$ of Assumption \ref{first} satisfy the parabolic Harnack inequality with a polynomial space-time scaling in a suitably consistent fashion, then the tightness condition of Assumption \ref{second} is immediately satisfied. The key result in proving that this is the case is provided by the following lemma, which demonstrates that the parabolic Harnack inequality implies the H\"{o}lder continuity of the transition density on graphs. The proof is an adaptation of \cite{BarHamllt}, Proposition 3.2, which deals with the case $\kappa=2$.

{\lem \label{philem} Fix $\kappa\geq 2$. Let $x\in V(G)$ and suppose that the parabolic Harnack inequality with constant $C_H$ holds for $Q_G(x,R,R^\kappa)$ for $R\geq s_G(x)$, where $s_G(x)$ is a positive integer depending on $x$. If $T^{1/\kappa}\geq 4R\geq 2s_G(x)$, then
\[\sup_{y\in B_G(x,R)}\left|q^G_{T}(x)-q_{T}^G(y)\right|\leq c \left(\frac{R}{T^{1/\kappa}}\right)^\theta\frac{1}{\nu^{G}(B_G(x,\tfrac{1}{4}T^{1/\kappa}))},\]
where $c$ and $\theta$ are constants depending only on $C_H$ taking values in $(0,\infty)$.}
\begin{proof} Let $x\in V(G)$ and suppose $T^{1/\kappa}\geq 4R\geq 2s_G(x)$. Set $T_0:=T+1$ and $R_0:=T_0^{1/\kappa}$. For $k\in \mathbb{N}$, define the space-time regions $Q(k):=[T_0-R_k^\kappa,T_0]\times B_G(x,R_k)$, and
\[Q_-(k):=[T_0-\tfrac{3}{4}R_k^\kappa,T_0-\tfrac{1}{2}R_k^\kappa]\times B_G(x,\tfrac{1}{2}R_k),\:Q_+(k):=[T_0-2^{-\kappa}R_k^\kappa,T_0)\times B_G(x,\tfrac{1}{2}R_k),\]
where $R_k:=2^{-k}R_0$. Since $\kappa\geq 2$, we have $T_0-2^{-\kappa}R_k^\kappa\geq T_0-\tfrac{1}{4}R_k^\kappa$. Consequently, if $R_k\geq s_G(x)$, then we can apply the parabolic Harnack inequality on $Q(k)$ to deduce that
\begin{equation}\label{phiq}
\sup_{(m,y)\in Q_-(k)}q^G_m(y)\leq C_H \inf_{(m,y)\in Q_+(k)} q^G_m(y).
\end{equation}
As in the proof of \cite{BarHamllt}, Proposition 3.2, it follows that we can bound $\mathrm{Osc}(q^G, Q_+(k))$ above by $(1-\tfrac{1}{2C_H})\mathrm{Osc}(q^G, Q(k))$, whenever $k$ satisfies $R_k\geq s_G(x)$. Noting that $Q_+(k)=Q(k+1)$, we can iterate this result to obtain that
\[\sup_{y\in B_G(x,R)}\left|q^G_{T}(x)-q_{T}^G(y)\right|\leq \left(1-\frac{1}{2C_H}\right)^{k-1}\mathrm{Osc}(q^G, Q(1)),\]
where $k$ is chosen to satisfy $R_k\geq 2R>R_{k+1}$. This implies that there exist constants $c_1$ and $\theta$, which depend only on $C_H$, that satisfy
\[\sup_{y\in B_G(x,R)}\left|q^G_{T}(x)-q_{T}^G(y)\right|\leq c_1 \left(\frac{R}{T^{1/\kappa}}\right)^\theta\sup_{(m,y)\in Q_-(1)}q^G_m(y).\]
Thus to complete the proof it suffices to bound the final term appropriately. Again applying (\ref{phiq}), we have that, for $m\in[T_0-\tfrac{3}{4}R_1^\kappa,T_0-\tfrac{1}{2}R_1^\kappa]$,
\[\nu^G(B_G(x,\tfrac{1}{4}R_0)) \sup_{y\in B_G(x,\tfrac{1}{4}R_0) }q^G_m(y)\leq C_H\int_{B_G(x,\tfrac{1}{4}R_0)}q^{G}_{T_0-2^{-\kappa}R_1^\kappa}(y)\nu^G(dy)\leq C_H.\]
The result follows.
\qed\end{proof}

To apply this H\"{o}lder continuity result, we make the following assumption on the graph sequence $(G^n)_{n\geq 1}$. Note that Assumption \ref{fourth}(b) is an extension of Assumption \ref{first}(b) and prevents elements of $V(G^n)$ being too far from $F$ for large $n$ (at least in bounded spatial regions).

\begin{assu}\label{fourth} In the setting of Assumption 1, suppose that the following statements are satisfied for some $\kappa\geq 2$ and $C_H<\infty$.
 \renewcommand{\labelenumi}{(\alph{enumi})}
\begin{enumerate}
\item Assumption \ref{first}(a) holds.
\item For every $r>0$,
\[\lim_{n\rightarrow\infty}V(G^n)\cap\overline{B}_E(\rho,r)= \overline{B}_F(\rho,r)\]
as $n\rightarrow \infty$ with respect to the usual Hausdorff topology on non-empty compact subsets of $(E,d_E)$.
\item For every $x\in V(G^n)$, $n\geq 1$, there exists a positive integer $s_{G^n}(x)$ such that the parabolic Harnack inequality with constant $C_H$ holds for $Q_{G^n}(x,R,R^\kappa)$ for $R\geq s_{G^n}(x)$. Moreover, suppose that there exists a dense subset $F^*$ of $F$ such that, for every $x\in F^*$, ${\alpha}(n)^{-1}s_{G^n}(g_n(x))\rightarrow0$.
\item As $n\rightarrow \infty$, we have $\alpha(n)^\kappa=O(\gamma(n))$.
\end{enumerate}
\end{assu}

The main result of this section is the following.

{\propn \label{42} If Assumption \ref{fourth} holds, then so does Assumption \ref{second}.}
\bigskip

Before we prove this result, however, we derive a lemma that describes a useful sequence of covers for balls of the form $B_{G^n}(\rho, \alpha(n)r)$.

{\lem \label{coverlem} Suppose Assumption \ref{fourth} holds. For every $r,\varepsilon>0$ there exists a finite set $\mathcal{X}\subseteq\overline{B}_F(\rho,r/c_1)\cap F^*$, where $c_1$ is the constant of the bound at (\ref{as1a}), and integer $n_0$ such that $(B_{G^n}(g_n(x),\alpha(n)\varepsilon))_{x\in\mathcal{X}}$ is a cover for $B_{G^n}(\rho,\alpha(n)r)$ whenever $n\geq n_0$.}
\begin{proof} Fix $r,\varepsilon>0$ and set $r_0:=r/c_1$, where $c_1$ is the constant of the bound at (\ref{as1a}). Choose $c_2$ and $n_0$ by Assumption \ref{first}(a) so that if $x,y\in V(G^n)\cap B_E(\rho,r_0+1)$ and $n\geq n_0$, then $d_{G^n}(x,y)\leq c_2\alpha(n)d_E(x,y)+\varepsilon \alpha(n)/4$. Furthermore, use Assumption \ref{fourth}(b) to find an integer $n_1\geq n_0$ such that
\begin{equation}\label{close1}
\sup_{x\in B_F(\rho,r_0+1)}d_E(x, g_n(x))<\varepsilon_0
\end{equation}
and
\begin{equation}\label{close2}
\sup_{x\in V(G^n)\cap\overline{B}_E(\rho,r_0)}d_E(x,\overline{B}_F(\rho,r_0))<\varepsilon_0
\end{equation}
for $n\geq n_1$, where $\varepsilon_0:=(\varepsilon/4c_2)\wedge 1$. As a final piece of information that we will need, note that, since $\overline{B}_F(\rho,r_0)$ is compact and $F^*$ is dense in $F$, there exists a finite collection $\mathcal{X}\subseteq \overline{B}_F(\rho,r_0)\cap F^*$ such that $(B_F(x,\varepsilon_0))_{x\in\mathcal{X}}$ is a cover for $\overline{B}_F(\rho,r_0)$. By (\ref{close1}), this implies that $(B_E(g_n(x),2\varepsilon_0))_{x\in\mathcal{X}}$ is a cover for $\overline{B}_F(\rho,r_0)$ for $n\geq n_1$.

Assume now that $n\geq n_1$ and let $x\in B_{G^n}(\rho,\alpha(n)r)$. Observe that Assumption \ref{first}(a) and (\ref{close2}) imply that $x\in B_E(\rho,r_0)$ and there exists a $y\in \overline{B}_F(\rho,r_0)$ such that $x\in B_E(y,\varepsilon_0)$ respectively. Thus, applying the final result of the previous paragraph, we have that $x\in B_E(g_n(y),3\varepsilon_0)$ for some $y\in\mathcal{X}$. Consequently, because $x,g_n(y)\in V(G^n)\cap B_E(\rho,r_0+1)$ and $d_E(x,g_n(y))<3\varepsilon_0$, it follows that $x\in B_{G^n}(g_n(y),\alpha(n)\varepsilon)$. Since the choice of $x\in B_{G^n}(\rho,\alpha(n)r)$ was arbitrary, the proof is complete.
\qed\end{proof}

\begin{proof}[Proof of Proposition \ref{42}] Fix $r, \varepsilon>0$, $I=[T_1,T_2]\subset(0,\infty)$ and suppose that $c_1$ is defined to be the constant of the bound at (\ref{as1a}). By Assumption \ref{fourth}(d), there exists a constant $c_2>0$ and an integer $n_0$ such that $\lfloor \gamma(n) T_1\rfloor^{1/\kappa} \geq c_2 \alpha(n)$ for $n\geq n_0$. Given these constants, Assumptions \ref{first}(a) and \ref{first}(b) imply that we can choose $r_0>0$ and an integer $n_1\geq n_0$ such that $V(G^n)\cap B_E(x,r_0)\subseteq B_{G^n}(g_n(x),\tfrac{1}{4}c_2\alpha(n))$ for every $x\in \overline{B}_F(\rho,r/c_1)$ and $n\geq n_1$. Furthermore, applying the compactness of $\overline{B}_F(\rho,r/c_1)$, we have that $c_3:=\inf_{x\in \overline{B}_F(\rho,r/c_1)}\nu(B_F(x,r_0))>0$. We use these constants to define
\[\delta:=\left(\frac{c_2^\theta c_3\varepsilon}{2^{\theta+2}c}\right)^{1/\theta}\wedge\left(\frac{c_2}{8}\right),\]
where $c$ and $\theta$ are the constants of Lemma \ref{philem} depending only on $C_H$.

By Lemma \ref{coverlem}, there exists a finite set $\mathcal{X}\subseteq \overline{B}_F(\rho,r/c_1)\cap F^*$ and an integer $n_2\geq n_1$ such that $(B_{G^n}(g_n(x),\alpha(n)\delta))_{x\in\mathcal{X}}$ is a cover for $B_{G^n}(\rho,\alpha(n)r)$ whenever $n\geq n_2$. Applying the finiteness of $\mathcal{X}$ and Assumptions \ref{first}(c) and \ref{fourth}(c), we are also able to deduce the existence of an integer $n_3\geq n_2$ such that $\max_{x\in\mathcal{X}}s_{G^n}(g_n(x))\leq 4\alpha(n) \delta$ and also
\[\beta(n)^{-1}\nu^n(B_E(x,r_0))\geq \tfrac{1}{2}\nu(B_E(x,r_0)),\hspace{20pt}\forall x\in\mathcal{X},\]
for $n\geq n_3$. In particular, if $n\geq n_3$, then we have $2s_{G^n}(g_n(x))\leq 8\alpha(n) \delta\leq \lfloor\gamma(n)t\rfloor^{1/\kappa}$ for every $t\in I$, and so we can apply Lemma \ref{philem} and our choice of constants to deduce that
\[\sup_{x\in\mathcal{X}}\sup_{y\in B_{G^n}(g_n(x),2\alpha(n)\delta)}\sup_{t\in I}\beta(n)\left|q^n_{\lfloor\gamma(n)t\rfloor}(g_n(x))-q_{\lfloor\gamma(n)t\rfloor}^n(y)\right|\leq \varepsilon /2\]
for every $n\geq n_3$. The proposition is a straightforward consequence of this inequality.
\qed\end{proof}

\section{Resistance estimates and tightness}\label{ressec}

As an alternative to the parabolic Harnack inequality, in this section we derive a sufficient condition for Assumption \ref{second} that involves an estimate of the resistance metric, which we now define. First, for a graph $G$, introduce an inner product $(\cdot,\cdot)_{G}$ on $\mathbb{R}^{V(G)}\times \mathbb{R}^{V(G)}$ by setting $(f,g)_{G}:=\sum_{x\in V(G)}f(x)g(x)\nu^G(\{x\})$. Use this and the discrete time generator $\mathcal{L}_G$ of the random walk $X^G$, as defined by (\ref{generator}), to construct a Dirichlet form $\mathcal{E}_G$ which satisfies
$\mathcal{E}_G(f,g):=-(\mathcal{L}_Gf,g)_{G}$. The domain of $\mathcal{E}_G$ is $\mathcal{F}_G:=\{f\in\mathbb{R}^{V(G)}:\:\mathcal{E}_G(f,f)<\infty\}$. For $x,y\in V(G)$, the resistance metric is defined by
\begin{equation}\label{resdef}
R_G(x,y):=\sup\left\{\frac{|f(x)-f(y)|^2}{\mathcal{E}_G(f,f)}:\:f\in\mathcal{F}_G,\:\mathcal{E}_G(f,f)>0\right\}.
\end{equation}
The following is proved as \cite{Barlow}, Proposition 4.25, see also \cite{Kigami}.

{\lem \label{res} The function $R_G$ is a metric on $V(G)$. Furthermore, for $f\in \mathcal{F}_G$,
\[(f(x)-f(y))^2\leq R_G(x,y)\mathcal{E}_G(f,f),\hspace{20pt}\forall x,y\in V(G).\]}
\bigskip

We will use this lemma to deduce oscillation bounds for $q_m^G$. To start with, observe that it is elementary to show that
\begin{equation}\label{equal}
\mathcal{E}_G(q_m^G,q_m^G)=q_{2m}^G(\rho)-q_{2m+2}^G(\rho)
\end{equation}
for every $m\geq 0$. This immediately implies that $\mathcal{E}_G(q_m^G,q_m^G)\leq q_{2m}^G(\rho)$. However, we will next prove a lemma demonstrating how to sharpen this bound. In the proof we use the notation $P_G$ to represent the linear operator defined from the transition probabilities $(P_G(x,y))_{x,y\in V(G)}$ of the simple random walk $X^G$ by
\[P_Gf(x)=\sum_{y\in V(G)}P_G(x,y)f(y),\]
for $f\in \mathbb{R}^{V(G)}$. Note that $P_G$ defines the usual random walk semigroup $(P_G^m)_{m\geq 0}$ and satisfies $P_G=\mathcal{L}_G+I_G$, where $I_G$ is the identity operator on $\mathbb{R}^{V(G)}$.

{\lem \label{energy} For every $m\geq 1$, we have $\mathcal{E}_G(q_{m}^G,q_{m}^G)\leq {2q_{2\lceil m/2\rceil}^G(\rho)}/{m}$.}
\begin{proof} Let us start by demonstrating that $(\mathcal{E}_G(q_m^G,q_m^G))_{m\geq 0}$ is a decreasing sequence. Applying (\ref{equal}) and the fact that $p^G_{m+n}(\rho,\rho)=(p_m^G(\rho,\cdot),p_n^G(\rho,\cdot))_G$, it is possible to deduce that
\[\mathcal{E}_G(q_m^G,q_m^G)-\mathcal{E}_G(q_{m+1}^G,q_{m+1}^G)=((I_G+P_G)(I_G-P_G^2)p^G_{m}(\rho,\cdot),(I_G-P_G^2)p^G_{m}(\rho,\cdot))_{G},\]
where we also apply the self-adjointness of $P_G$ with respect to $(\cdot,\cdot)_G$. Since $P_G$ is stochastic, it can easily be checked that $((I_G+P_G)f,f)_G\geq 0$ for every $f\in\mathbb{R}^{V(G)}$, and therefore $\mathcal{E}_G(q_m^G,q_m^G)\geq\mathcal{E}_G(q_{m+1}^G,q_{m+1}^G)$, as desired.

Again applying (\ref{equal}) we see that $\sum_{m=M}^{2M-1}\mathcal{E}_G(q_m^G,q_m^G)=q_{2M}^G(\rho)-q_{4M}^G(\rho)\leq q_{2M}^G(\rho)$. Since the summands of the left-hand side are decreasing in $m$, we consequently have that $M\mathcal{E}_G(q^G_{2M-1},q^G_{2M-1})\leq q_{2M}^G(\rho)$, and the result follows from this.
\qed\end{proof}

We now describe how to bound $q_{2m}^G(\rho)$ in terms of the volume growth about $\rho$ of the graph $G$ with respect to the resistance metric. Define a function $V_G:\mathbb{R}_+\rightarrow\mathbb{R}_+$ by setting
\[V_G(r):=\nu^G(\{x:\:R_G(\rho,x)\leq r\}),\]
so that $V_G(r)$ represents the volume of the closed ball around $\rho$ of radius $r$ with respect to the resistance metric. Set $h_G(r):=rV_G(r)$, and define the  right-continuous inverse of $h_G$ by
\begin{equation}\label{hinv}
h_G^{-1}(m):=\sup\{r:\:h_G(r)\leq m\}.
\end{equation}
The proof of the following result is a simple adaptation of \cite{BCK}, Proposition 3.2. We do, however, continue to include the proof in order to demonstrate the universality of the constant in the resultant upper bound.

{\lem For every $m\geq 1$, we have $q_{2m}^G(\rho)\leq{3h_G^{-1}(m)}/{m}$.}
\begin{proof} Since, for $r>0$, $\sum_{{x\in V(G):R_G(\rho,x)\leq r}}q_{2m}^G(x)\nu^G(\{x\})\leq 1$, there must exist an $x\in V(G)$ such that $R_G(\rho,x)\leq r$ and also $q^G_{2m}(x)\leq V_G(r)^{-1}$. Hence we obtain
\begin{eqnarray*}
q^G_{2m}(\rho)^2&\leq & 2 q_{2m}^G(x)^2+2(q_{2m}^G(x)-q_{2m}^G(\rho))^2\\
&\leq& 2 V_G(r)^{-2}+2r\mathcal{E}_G(q_{2m}^G,q_{2m}^G)\\
&\leq & 2 V_G(r)^{-2}+\frac{2r q_{2m}^G(\rho)}{m},
\end{eqnarray*}
where we have applied Lemma \ref{res} for the second inequality and Lemma \ref{energy} for the third. This quadratic inequality implies that
\[q^G_{2m}(\rho)\leq \frac{r}{m}+\frac{1}{2}\sqrt{\frac{4r^2}{m^2}+\frac{4}{V_G(r)^2}}
\leq \frac{2r}{m}+\frac{1}{V_G(r)}.\]
The result follows on choosing $r=h^{-1}_G(m)$.
\qed\end{proof}

Combining the three previous lemmas we obtain the following result.

{\propn \label{resprop} For $m\geq 1$,
\[(q^G_m(x)-q^G_m(y))^2\leq \frac{12R_G(x,y)h_G^{-1}(\lceil m/2\rceil)}{m^2},\hspace{20pt}\forall x,y\in V(G).\]}
\bigskip

Application of the above bound relies on being able to adequately control the resistance between points in $V(G)$ and the volume growth with respect to the resistance metric, which is not always possible. However, as we shall demonstrate in Section \ref{examples}, there are classes of graphs for which we can make use of this result, most notable amongst these are nested fractal graphs and graph trees. More specifically, in the case when the resistance metric $R_{G^n}$ is bounded above by a power of the shortest path metric $d_{G^n}$ for graphs in the sequence $(G^n)_{n\geq 1}$, the tightness condition of Assumption \ref{second} follows directly from Assumptions \ref{first}(a) and \ref{first}(c) whenever the space, volume and time scaling factors are related in a way we now describe.

\begin{assu} \label{fifth} Suppose that, in the setting of Assumption \ref{first}, $(G^n)_{n\geq 1}$ is a sequence of graphs for which Assumptions \ref{first}(a) and \ref{first}(c) hold and, for some $\kappa\in(0,\infty)$, there exist constants $c_1,c_2,c_3\in(0,\infty)$ and an integer $n_0$ such that
\[R_{G^n}(x,y)\leq c_1d_{G^n}(x,y)^\kappa,\hspace{20pt}\forall x,y\in V(G^n),\]
and also
\begin{equation}\label{er}
c_2\gamma(n)\leq \alpha(n)^\kappa\beta(n)\leq c_3\gamma(n),
\end{equation}
for $n\geq n_0$.
\end{assu}
\bigskip

Under this assumption we can bound functions of the form $\alpha(n)^{-\kappa}h_{G^n}^{-1}(\gamma(n)\cdot)$, $n\geq 1$,  uniformly over compact intervals. In the proof of the following result, we consider the function $v:\mathbb{R}_+\rightarrow\mathbb{R}_+$ which satisfies $v(r):=\nu(B_E(\rho,r))$.

{\lem\label{hbound} Suppose Assumption \ref{fifth} holds. For any compact interval $I\subset (0,\infty)$,
\[\limsup_{n\rightarrow \infty}\sup_{t\in I}\alpha(n)^{-\kappa}h_{G^n}^{-1}(\gamma(n)t)<\infty.\]}
\begin{proof} Fix $I=[T_1,T_2]\subset(0,\infty)$. Let $c_1, c_2$ and $n_0$ be chosen such that $c_1d_{G^n}(x,y)^\kappa\geq R_{G^n}(x,y)$ for every $x,y\in V(G^n)$ and $\alpha(n)^\kappa\beta(n)\geq c_2\gamma(n)$ for $n\geq n_0$, which is possible by Assumption \ref{fifth}. Define $R\in (0,\infty)$ to be a constant satisfying $c_1c_2R^\kappa v(R)>T_2$. By Assumption \ref{first}(a), there also exists a constant $c_3\geq 1$ and integer $n_1\geq n_0$ such that $d_{G^n}(x,y)\leq c_3\alpha(n)d_E(x,y)+\tilde{\alpha}(n)$ for every $x,y\in V(G^n)\cap B_E(\rho,R)$, $n\geq n_1$. Furthermore, define $n_2\geq n_1$ to be an integer such that $\tilde{\alpha}(n)\leq c_3\alpha(n)R$ for $n\geq n_2$, then we have that
\begin{equation}\label{subset}
V(G^n)\cap B_E(\rho,R)\subseteq B_{G^n}(\rho,2c_3\alpha(n)R)\subseteq\left\{x\in V(G^n):\:R_{G^n}(\rho,x)\leq c_4 \alpha(n)^\kappa\right\},
\end{equation}
for $n\geq n_2$, where $c_4:=c_1(2c_3R)^\kappa$.

Now, by Assumption \ref{first}(c), there exists an $n_3\geq n_2$ such that $\beta(n)^{-1}\nu^n(B_E(\rho,R))$ is bounded below by $2^{-\kappa}v(R)$ for every $n\geq n_3$. Consequently, applying (\ref{subset}), we have that $V_{G^n}(c_4\alpha(n)^\kappa)> 2^{-\kappa}\beta(n)v(R)$, for every $n\geq n_3$. It immediately follows that $h_{G^n}(c_4\alpha(n)^\kappa)> c_1 c_2c_3^\kappa\gamma(n)R^\kappa v(R)$ for $n\geq n_3$. Thus
\[\sup_{t\in I}h_{G^n}^{-1}(\gamma(n)t)\leq h_{G^n}^{-1}(\gamma(n)T_2)\leq  h_{G^n}^{-1}(c_1c_2c_3^\kappa\gamma(n)R^\kappa v(R))\leq c_4\alpha(n)^\kappa,\]
for $n\geq n_3$, which completes the proof.
\qed\end{proof}

We now arrive at the first main result of this section.

{\propn \label{res2} If Assumption \ref{fifth} holds, then so does Assumption \ref{second}.}
\begin{proof} Fix an interval $I=[T_1,T_2]\subset(0,\infty)$ and $r>0$. It is straightforward to obtain from Proposition \ref{resprop} and Lemma \ref{hbound} the existence of a finite constant $c_1$ and an integer $n_0$ such that
\[\sup_{\substack{x,y\in V(G^n):\\d_{G^n}(x,y)\leq \alpha(n)\delta}}\sup_{t\in I}\beta(n)^2\left|q_{\lfloor\gamma(n)t\rfloor}^n(x)-q_{\lfloor\gamma(n)t\rfloor}^n(y)\right|^2\leq c_1\left(\frac{\alpha(n)^\kappa\beta(n)}{\gamma(n)}\right)^2\delta^\kappa,\]
for every $\delta>0$ and $n\geq n_0$. Hence the inequality at (\ref{er}) implies the proposition.
\qed\end{proof}

To complete this section, let us remark that the bounds at (\ref{er}) can be interpreted in terms of the random walk version of the Einstein relation, which explains how the time, resistance and volume scaling exponents for random walks on graphs are related (see \cite{Telcs} for background). In particular, if we assume that for a graph $G$ the resistance satisfies $R_G\asymp d_G^\kappa $ (where $\asymp$ is taken to mean ``bounded above and below by constant multiples of'') and the volume satisfies $\nu^G(B_G(x,r))\asymp r^d$, then it is possible to deduce that
\[\mathbf{E}_\rho^GT_G(\rho,r)\asymp r^{\kappa+d},\]
where $T_G(\rho,r):=\min\{m\geq 0:\:X^G_m\not\in B_G(\rho,r)\}$ is the exit time of $X^G$ from $B_G(\rho,r)$, and $\mathbf{E}_\rho^G$ is the expectation under $\mathbf{P}_\rho^G$ (see \cite{BCK}, Proposition 3.4, for example). Thus, if such polynomial relations hold uniformly for the graphs in the sequence $(G^n)_{n\geq 1}$, then, from the scaling considerations of Assumption \ref{first}, one might expect to be able to conclude that
\[\gamma(n)\asymp \mathbf{E}_\rho^{G^n}T_{G^n}(\rho,\alpha(n))\asymp \alpha(n)^{\kappa+d},\hspace{20pt}\beta(n)\asymp\nu^{n}(B_{G^n}(x,\alpha(n)))\asymp \alpha(n)^d,\]
which would imply that $\gamma(n)\asymp \alpha(n)^\kappa\beta(n)$, as required for (\ref{er}) to hold.

\section{Two-spatial parameter local limit theorems}\label{2param}

So far we have considered the asymptotics of the transition densities of the simple random walks on graphs in a sequence $(G^n)_{n\geq 0}$ when the relevant processes are started from a fixed point $\rho$. We now provide conditions that will allow us to extend these results uniformly to arbitrary starting points and deduce local limit theorems for the two-spatial parameter functions $(q^n_m(x,y))_{x,y\in V(G^n), m\geq 0}$, $n\geq 1$, as defined at (\ref{qdef}). In this section, we assume that $X=((X_t)_{t\geq 0}, \mathbf{P}_x,x\in F)$ is a conservative $\nu$-symmetric Markov diffusion on $F$, with a transition density $(p_t(x,y))_{x,y\in F, t>0}$ which is jointly continuous in $(t,x,y)$. The extensions of Assumptions \ref{first} and \ref{second} we apply are the following.

\begin{assu1hat}{\it In the setting of Assumption \ref{first}, suppose that Assumptions \ref{first}(a), \ref{first}(b) and \ref{first}(c) are satisfied. Moreover, suppose that there exists a dense subset $F^*$ of $F$ such that, for any compact interval $I\subset (0,\infty)$, $x\in F^*$, $y\in F$, and $r>0$,
\begin{equation}\label{2p}\lim_{n\rightarrow\infty}\mathbf{P}_{g_n(x)}^{G^n}\left(X^n_{\lfloor \gamma(n)t\rfloor}\in B_E(y,r)\right)=\mathbf{P}_{x}\left(X_t\in B_E(y,r)\right)\end{equation}
uniformly for $t\in I$.}
\end{assu1hat}

\begin{assu2hat}{\it In the setting of Assumption \ref{first}, suppose that, for any compact interval $I\subset (0,\infty)$ and $r>0$,
\[\lim_{\delta\rightarrow 0}\limsup_{n\rightarrow\infty}\sup_{\substack{x,y,z\in B_{G^n}(\rho,\alpha(n)r):\\d_{G^n}(y,z)\leq \alpha(n)\delta}}\sup_{t \in I}\beta(n)\left|q_{\lfloor\gamma(n)t\rfloor}^{n}(x,y)-q_{\lfloor\gamma(n)t\rfloor}^{n}(x,z)\right|=0.\]}
\end{assu2hat}

We now prove our main two-spatial parameter local limit theorem, which is a variation of Theorem \ref{bdd}.

{\thm \label{bdd2} Fix a compact interval $I\subset (0,\infty)$ and $r>0$. Suppose Assumptions \^{1} and \^{2} hold, then
\[\lim_{n\rightarrow\infty}\sup_{x,y\in B_F(\rho,r)}\sup_{t\in I}\left|\beta(n)q^n_{\lfloor\gamma(n)t\rfloor}(g_n(x),g_n(y))-p_t(x,y)\right|=0.\]}
\begin{proof} Fix a compact interval $I\subset (0,\infty)$. For any $x\in F^*$, $y\in F$, we can prove that
\begin{equation}\label{point2}
\lim_{n\rightarrow\infty}\sup_{t\in I}\left|\beta(n)q^n_{\lfloor\gamma(n)t\rfloor}(g_n(x),g_n(y))-p_t(x,y)\right|=0
\end{equation}
by following a proof similar to the proof of Proposition \ref{point}. To extend this result to hold uniformly over $x,y\in B_F(\rho,r)$, we can proceed as in the proof of Theorem \ref{bdd} by first applying (\ref{point2}) to deduce the result holds uniformly over a suitably chosen finite set $\mathcal{X}\subseteq \overline{B}_F(\rho,r)$, and then using the tightness condition of Assumption \^{2} and the continuity of $(p_t(x,y))$ to extend this to the whole of $B_F(\rho,r)$. Note that in order to strengthen (\ref{point2}) in this way, one should choose $\mathcal{X}\subseteq\overline{B}_F(\rho,r)\cap F^*$, which is possible by the denseness of $F^*$ in $F$.
\qed\end{proof}

If we suppose that $(F,d_F)$ satisfies the midpoint property and the obvious extensions to the transition density decay conditions of Assumption \ref{third} hold, then, by following a proof similar to that of Theorem \ref{unbdd}, it is possible to extend this result to demonstrate that $\beta(n)q^n_{\lfloor\gamma(n)t\rfloor}(g_n(x),g_n(y))$ converges uniformly to $p_t(x,y)$ over $(t,x,y)\in [T_1,\infty)\times B_F(\rho,R)\times F$, for any $T_1, R>0$. To prove uniform convergence of the transition densities on $[T_1,\infty)\times F^2$ in general, however, seems to require some uniform control over space of the convergence of measures and processes of Assumption \ref{first}(c) and Assumption \^{1}.

We now extend Propositions \ref{42} and \ref{res2} to show that the parabolic Harnack inequality of Assumption \ref{fourth} and the resistance estimates of Assumption \ref{fifth} imply the uniform tightness condition of Assumption \^{2}.

{\propn If Assumption \ref{fourth} holds, then Assumption \^{2} holds.}
\begin{proof} Observe that the proof of Lemma \ref{philem} only depended on $(p^G_m(\rho,y))$ being a caloric function of $(m,y)$ and the fact that $\int_{V(G)}q^G_m(y)\nu^G(dy)\leq1$ for any $m\geq 1$. Hence, because the same is true of $(p^G_m(x,y))$ and $(q^G_m(x,y))$ for any $x\in V(G)$, we are able to deduce that if $T^{1/\kappa}\geq 4R\geq 2s_G(y)$, then
\[\sup_{x\in V(G)}\sup_{z\in B_G(y,R)}\left|q^G_{T}(x,y)-q_{T}^G(x,z)\right|\leq c \left(\frac{R}{T^{1/\kappa}}\right)^\theta\frac{1}{\nu^{G}(B_G(y,\tfrac{1}{4}T^{1/\kappa}))},\]
where $c,\theta\in (0,\infty)$ are the constants of Lemma \ref{philem} depending only on $C_H$. Consequently, if we fix $r,\varepsilon>0$, a compact interval $I\subset(0,\infty)$, and then choose $\delta>0$ and $\mathcal{X}\subseteq F^*$ as in the proof of Proposition \ref{42}, we obtain the existence of an integer $n_0$ such that
\[\sup_{x\in V(G^n)}\sup_{y\in\mathcal{X}}\sup_{z\in B_{G^n}(g_n(y),2\alpha(n)\delta)}\sup_{t\in I}\beta(n)\left|q^n_{\lfloor\gamma(n)t\rfloor}(x,g_n(y))-q_{\lfloor\gamma(n)t\rfloor}^n(x,z)\right|\leq \varepsilon /2\]
and $(B_{G^n}(g_n(y), \alpha(n)\delta))_{y\in \mathcal{X}}$ is a cover for $B_{G^n}(\rho,\alpha(n)r)$, whenever $n\geq n_0$. The proposition follows.
\qed\end{proof}

{\propn \label{52hat} If Assumption \ref{fifth} holds, then Assumption \^{2} holds}
\begin{proof} Generalising the notation from Section \ref{ressec}, let $V_G(x,r):=\nu^G(\{y:\:R_G(x,y)\leq r\})$, $h_G(x,r):=rV_G(x,r)$, and define $h_G^{-1}(x,m)$ by a formula analogous to (\ref{hinv}). Fix a compact interval $I\subset (0,\infty)$. By applying the estimate of Proposition \ref{resprop} to $q_m^n(x,\cdot)$ instead of $q_m^n(\cdot)$ and the bounds of Assumption \ref{fifth}, we can deduce that there exists finite constants $c_1$ and $c_2$ and integer $n_0$ such that
\begin{eqnarray*}
\lefteqn{\sup_{\substack{x,y,z\in B_{G^n}(\rho,\alpha(n)r):\\d_{G^n}(y,z)\leq \alpha(n)\delta}}\sup_{t \in I}\beta(n)^2\left|q_{\lfloor\gamma(n)t\rfloor}^{n}(x,y)-q_{\lfloor\gamma(n)t\rfloor}^{n}(x,z)\right|^2}\\
&\hspace{80pt}\leq& c_1\delta^\kappa \sup_{x\in B_{G^n}(\rho,\alpha(n)r)}\alpha(n)^{-\kappa}h^{-1}_{G^n}(x,c_2\gamma(n)),
\end{eqnarray*}
for every $\delta>0$ and $n\geq n_0$. Thus to complete the proof it suffices to obtain an asymptotic bound for the supremum in this expression, which can be achieved by a simple extension of Lemma \ref{hbound}.
\qed\end{proof}

\section{Local limit theorems for random weights}\label{randsec}

We now explain how Theorems \ref{bdd} and \ref{unbdd} can be generalised to the case where the weight functions on the graphs in the sequence $(G^n)_{n\geq 1}$ are chosen randomly from a law $\mathbb{P}_\mu$, a probability measure on $(0,\infty)^{\cup_{n\geq 0}E(G^n)}$. The adaptations of Assumptions 1 and 2 that we will apply are the following probabilistic versions.

\begin{assu1r}{\it In the setting of Assumption \ref{first}, suppose that Assumptions \ref{first}(a) and \ref{first}(b) hold. Moreover, suppose that, for every $x\in F$, $r,\varepsilon>0$,
\begin{equation}\label{randommeasconv}\lim_{n\rightarrow\infty}\mathbb{P}_{\mu}\left(\left|\beta(n)^{-1}\nu^{n}(B_E(x,r))-\nu(B_E(x,r))\right|>\varepsilon\right)=0,
\end{equation}
and, for any compact interval $I\subset (0,\infty)$, $x\in F$ and $r, \varepsilon>0$,
\begin{equation}\label{randomprocconv}
\lim_{n\rightarrow\infty}\mathbb{P}_{\mu}\left(\sup_{t\in I}\left|\mathbf{P}_{\rho}^{G^n}\left(X^n_{\lfloor \gamma(n)t\rfloor}\in B_E(x,r)\right)-\int_{B_E(x,r)}q_t(y)\nu(dy)\right|>\varepsilon\right)=0.
\end{equation}}
\end{assu1r}

\begin{assu2r}{\it In the setting of Assumption \ref{first}, suppose that, for any compact interval $I\subset (0,\infty)$ and $r,\varepsilon>0$,
\[\lim_{\delta\rightarrow 0}\limsup_{n\rightarrow\infty}\mathbb{P}_{\mu}\left(
\sup_{\substack{x,y\in B_{G^n}(\rho,\alpha(n)r):\\d_{G^n}(x,y)\leq \alpha(n)\delta}}\sup_{t \in I}\beta(n)\left|q_{\lfloor\gamma(n)t\rfloor}^{n}(x)-q_{\lfloor\gamma(n)t\rfloor}^{n}(y)\right|>\varepsilon\right)=0.\]}
\end{assu2r}

These assumptions allow us to prove the subsequent probabilistic local limit theorem.

{\thm Fix a compact interval $I\subset (0,\infty)$ and $r,\varepsilon>0$. Suppose Assumptions 1R and 2R hold, then
\[\lim_{n\rightarrow\infty}\mathbb{P}_{\mu}\left(\sup_{x\in B_F(\rho,r)}\sup_{t\in I}\left|\beta(n)q^n_{\lfloor\gamma(n)t\rfloor}(g_n(x))-q_t(x)\right|>\varepsilon\right)=0.\]}
\begin{proof} Fix a compact interval $I\subset (0,\infty)$, $x\in F$, $\varepsilon,\eta>0$, set $r=d_E(\rho,x)$, and choose $c$ and $n_0$ as in the proof of Proposition \ref{point}. By Assumption 2R, there exists an $r_0$ small enough and integer $n_1\geq n_0$ such that
\begin{equation}\label{bound}
\mathbb{P}_{\mu}\left(
\sup_{\substack{x,y\in B_{G^n}(\rho,c\alpha(n)(r+2)):\\d_{G^n}(x,y)\leq 3c\alpha(n)r_0}}\sup_{t \in I}\beta(n)\left|q_{\lfloor\gamma(n)t\rfloor}^{n}(x)-q_{\lfloor\gamma(n)t\rfloor}^{n}(y)\right|>\varepsilon/2\right)<\eta,
\end{equation}
for every $n\geq n_1$, and (\ref{hkcont}) holds. Consider $J(t,n)$ and $J_i(t,n)$, $i=1,\dots,4$ as in the proof of Proposition \ref{point}, and note that
\[\mathbb{P}_\mu\left(\sup_{t\in I}|J_2(t,n)|>4\varepsilon'\right)\leq \mathbb{P}_\mu\left(\sup_{t\in I}|J(t,n)|>\varepsilon'\right)+\sum_{i=1,3,4}\mathbb{P}_\mu\left(\sup_{t\in I}|J_i(t,n)|>\varepsilon'\right),\]
where $\varepsilon':=\varepsilon\nu(B_E(x,r_0))$. By (\ref{hkcont}), the term involving $J_4$ is equal to 0. Furthermore, by Assumption 1R, we can also choose $n_2\geq n_1$ large enough so that the terms featuring $J$ and $J_3$ are bounded above by $\eta$ for $n\geq n_2$. Note that to extend the convergence at (\ref{randomprocconv}) to $\lim_{n\rightarrow\infty}\mathbb{P}_\mu\left(\sup_{t\in I}|J(t,n)|>\varepsilon\right)=0$ for any $\varepsilon>0$, we apply a simple adaptation of the argument appearing in the proof of Lemma \ref{gen}. For the $J_1$ term, we first apply the upper bound for $J_1$ appearing at (\ref{j1bound}) and then (\ref{bound}), to deduce that
\[\mathbb{P}_\mu\left(\sup_{t\in I}|J_1(t,n)|>\varepsilon'\right)\leq \eta + \mathbb{P}_\mu\left(\beta(n)^{-1}\nu^n(B_E(x,r_0))> 2 \nu(B_E(x,r_0))\right),\]
for $n\geq n_2$. Thus Assumption 1R implies the existence of an integer $n_3\geq n_2$ such that $\mathbb{P}_\mu\left(\sup_{t\in I}|J_1(t,n)|>\varepsilon'\right)\leq 2\eta$ for $n\geq n_3$. Consequently, for some integer $n_4\geq n_3$, we have that
$\mathbb{P}_{\mu}(\sup_{t\in I}|\beta(n)q^n_{\lfloor\gamma(n)t\rfloor}(g_n(x))-q_t(x)|>8\varepsilon) \leq 6\eta$, for $n\geq n_4$.

We now explain how to generalise this point-wise result to hold uniformly over balls in $F$. Fix $r,\varepsilon,\eta>0$ and apply Assumption 2R to choose $r_0$ as above, so that (\ref{hkcont}) and (\ref{bound}) both hold for large $n$. As in the proof of Theorem \ref{bdd}, let $\mathcal{X}\subseteq \overline{B}_F(\rho,r)$ be a finite set such that $(B_F(x,r_0))_{x\in \mathcal{X}}$ is a cover for ${B}_F(\rho,r)$, then we can apply the first part of the proof to deduce that
\[\mathbb{P}_{\mu}\left(\sup_{x\in \mathcal{X}}\sup_{t\in I}\left|\beta(n)q^n_{\lfloor\gamma(n)t\rfloor}(g_n(x))-q_t(x)\right|>\varepsilon\right)<\eta,\]
for large enough $n$. The theorem follows from this by applying the continuity and tightness results of (\ref{hkcont}) and (\ref{bound}) respectively, similarly to the proof of Theorem \ref{bdd}.
\qed\end{proof}

The decay condition of Assumption 3 has the following analogous formulation.

\begin{assu3r}{\it The metric space $(F,d_F)$ has the midpoint property. Furthermore, the following conditions are fulfilled.
 \renewcommand{\labelenumi}{(\alph{enumi})}
\begin{enumerate}
\item Assumption \ref{third}(a) holds.
\item In the setting of Assumption \ref{first}, for every $\varepsilon>0$,
\[\lim_{t\rightarrow\infty}\limsup_{n\rightarrow\infty}\mathbb{P}_\mu\left(\sup_{x\in V(G^n)} \beta(n)q_{\lfloor \gamma(n)t\rfloor}^n(x)>\varepsilon\right)=0,\]
and, for any compact interval $I\subset (0,\infty)$, $\varepsilon>0$,
\[\lim_{r\rightarrow\infty}\limsup_{n\rightarrow\infty}\mathbb{P}_\mu\left(\sup_{x\in V(G^n)\backslash B_{G^n}(\rho,\alpha(n)r)} \sup_{t\in I}\beta(n)q_{\lfloor \gamma(n)t\rfloor}^n(x)>\varepsilon\right)=0.\]
\end{enumerate}}
\end{assu3r}

If this assumption is satisfied, then by decomposing time and space as in the proof of Theorem \ref{unbdd} we are able to deduce the corresponding result for random weights. Since the proof is a straightforward adaptation of the proof of Theorem \ref{unbdd}, we omit it.

{\thm Fix $T_1>0$. Suppose Assumptions 1R, 2R and 3R hold, then, for every $\varepsilon>0$,
\[\lim_{n\rightarrow\infty}\mathbb{P}_{\mu}\left(\sup_{x\in F}\sup_{t\geq T_1}\left|\beta(n)q^n_{\lfloor\gamma(n)t\rfloor}(g_n(x))-q_t(x)\right|>\varepsilon\right)=0.\]}
\bigskip

We can also extend the two-spatial parameter local limit theorems of Section \ref{2param} to random weights; see Section \ref{nf} for such a result.

\section{Examples}\label{examples}

To demonstrate the applicability of our local limit theorems, in this section we present a range of examples for which we can check that our assumptions hold.

\subsection{Lattice graphs}

In \cite{BarHamllt} local limit theorems were proved for an infinite subgraph $\mathcal{G}$ of the integer lattice $\mathbb{Z}^d$ fulfilling certain conditions, including a version of the parabolic Harnack inequality related to our Assumption \ref{fourth} with $\kappa=2$. It is easy to check that if we set $G^n:=n^{-1/2}\mathcal{G}$ for a graph $\mathcal{G}$ satisfying Assumption 4.4. of \cite{BarHamllt},
by which we mean that
\[V(G^n)=n^{-1/2}V(\mathcal{G}),\hspace{20pt}E(G^n):=\{\{n^{-1/2}x,n^{-1/2}y\}:\:\{x,y\}\in E(\mathcal{G})\},\]
\[\mu_{xy}^{G^n}=\mu^\mathcal{G}_{(n^{1/2}x)(n^{1/2}y)},\]
and define $\rho(G^n)=0$, then our Assumptions \ref{first}, \ref{third} and \ref{fourth} hold with: $(E,d_E)=(F,d_F)=(\mathbb{R}^d,|\cdot-\cdot|_\infty)$, where $|\cdot|_\infty$ is the usual $L^\infty$ norm in $\mathbb{R}^d$; $\rho=0$; $\nu$ equal to Lebesgue measure on $\mathbb{R}^d$; for some constant $c_1\in(0,\infty)$,
\[q_t(x)=\frac{1}{(2\pi c_1 t)^{d/2}}e^{-|x|^2/2c_1t};\]
and
\[\alpha(n)=n^{1/2},\hspace{20pt}\beta(n)=c_2n^{d/2},\hspace{20pt}\gamma(n)=n,\]
for some constant $c_2\in(0,\infty)$. Since, by Proposition \ref{42}, Assumption \ref{fourth} implies Assumption \ref{second}, we can apply Theorem \ref{unbdd} to verify the local limit theorem proved as \cite{BarHamllt}, Theorem 4.5. Examples of $\mathcal{G}$ to which such an argument applies include: the unweighted ($\mu^{\mathcal{G}}_{xy}=1$ for $\{x,y\}\in E(\mathcal{G})$) integer lattice $\mathbb{Z}^d$; typical supercritical percolation clusters; and typical realisations of the weighted graph generated by the random conductance model on $\mathbb{Z}^d$ in the case when the conductances are uniformly bounded away from 0 and $\infty$. See \cite{BarHamllt} for details.

\subsection{Graph trees converging to the continuum random tree}

To describe a scaling limit result for ordered graph trees, we will use the now well-established connection between trees and excursions (see \cite{Aldous2}, \cite{rrt}, for example).  First, let $(\mathcal{T}_n)_{n\geq 2}$ be a collection of (rooted) ordered graph trees, where $\mathcal{T}_n$ has $n$ vertices. For each $n$, define the function $\hat{w}_n:\{1,\dots,2n-1\}\rightarrow\mathcal{T}_n$ to be the depth-first search around $\mathcal{T}_n$ (see \cite{Aldous3} for a definition). Extend $\hat{w}_n$ so that $\hat{w}_n(0)=\hat{w}_n(2n)=\rho_n$, where $\rho_n$ is the root of $\mathcal{T}_n$. Define the search-depth process $w_n$ by $w_n(i/2n):=d_{\mathcal{T}_n}(\rho_n,\hat{w}_n(i))$ for $0\leq i\leq 2n$, where $d_{\mathcal{T}_n}$ is the graph distance on $\mathcal{T}_n$. Also, extend the definition of $w_n$ to the whole of the interval $[0,1]$ by linear interpolation, so that $w_n$ takes values in $C([0,1],\mathbb{R}_+)$.

Now, set $\mathcal{W}:=\{w\in C([0,1],\mathbb{R}_+):\:w(x)=0\Leftrightarrow x\in\{0,1\}\}$. For each $w\in \mathcal{W}$, construct on $[0,1]$ a distance by setting $d_w(s,t):=w(s)+w(t)-2m_w(s,t)$, where $m_w(s,t):=\inf\{w(r):\:r\in[s\wedge t,s\vee t]\}$, and an equivalence relation by supposing $s\sim t$ if and only if $d_w(s,t)=0$. If $\mathcal{T}_w:=[0,1]/\sim$ and $d_{\mathcal{T}_w}([s],[t]):=d_w(s,t)$, where $[s]$ is the equivalence class containing $s$, then it is possible to check that $(\mathcal{T}_w,d_{\mathcal{T}_w})$ is a compact real tree (see \cite{LegallDuquesne} for a definition of a real tree and proof of this result). The root $\rho_w$ of the tree $\mathcal{T}_w$ is defined to be the equivalence class $[0]$. Furthermore, if $\nu_w(A)$ is defined to be the standard one-dimensional Lebesgue measure of the set $\{s\in[0,1]:\:[s]\in A\}$ for Borel measurable $A\subseteq \mathcal{T}_w$, then $\nu_w$ is a Borel probability measure on $(\mathcal{T}_w,d_{\mathcal{T}_w})$. Since the support of $\nu_w$ is the whole of $\mathcal{T}_w$, it is possible to apply \cite{Kigamidendrite}, Theorem 5.4 to deduce the existence of a reversible strong Markov diffusion on $\mathcal{T}_w$, $X^w$ say, which has $\nu_w$ as its invariant measure. Moreover, by applying the argument of \cite{Croydoncrt}, Section 8, we can suppose that $X^w$ is Brownian motion on $(\mathcal{T}_w,d_{\mathcal{T}_w},\nu_w)$, as defined in Section 5 of \cite{Aldous2}.

The continuum random tree is the random compact real tree $\mathcal{T}_W$ that results when $W$ is the Brownian excursion, normalised to have length one. If the sequence $(w_n)_{n\geq 2}$ converges to a typical realisation of $W$, $w$ say, in $C([0,1],\mathbb{R}_+)$, then it was shown in \cite{Croydoncbp} that it is also possible to describe the scaling limits of the vertex sets $V(\mathcal{T}_n)$, the measures $\nu^{n}:=\nu^{\mathcal{T}_n}$ and the discrete time simple random walks $X^n:=X^{\mathcal{T}_n}$ in terms of the corresponding continuum objects $\mathcal{T}_w$, $\nu_w$ and $X^w$ by embedding in to a common metric space. In the following result, which is a minor restatement of \cite{Croydoncbp}, Theorem 1.1, the space $l^1$ is the Banach space of infinite sequences of real numbers equipped with the norm $\|x\|:=\sum_{i=1}^\infty |x_i|$.

{\propn \label{prop}There exists a set $\mathcal{W}^*\subseteq C([0,1],\mathbb{R}_+)$ such that $W\in \mathcal{W}^*$ almost-surely, and if $n^{-1/2}w_n\rightarrow w$ in $C([0,1],\mathbb{R}_+)$ for some $w\in\mathcal{W}^*$, then there exists, for each $n$, an isometric embedding $\phi_n$ of $(V(\mathcal{T}_n),d_{\mathcal{T}_n})$ into $l^1$ and also an isometric embedding $\phi$ of $(\mathcal{T}_w,d_{\mathcal{T}_w})$ into $l^1$ such that:
\begin{itemize}
  \item $\phi_n(\rho_n)=0$, for every $n\geq 2$, and also $\phi(\rho_w)=0$.
  \item $n^{-1/2}\phi_n(V(\mathcal{T}_n))\rightarrow \phi(\mathcal{T}_w)$ with respect to the Hausdorff topology on compact subsets of $l^1$.
  \item $(2n)^{-1}\nu^n(\phi_n^{-1}(n^{1/2}\cdot))\rightarrow \nu_w(\phi^{-1}(\cdot))$ weakly as Borel probability measures on $l^1$.
  \item $\left(n^{-1/2}\phi_n(X^n_{\lfloor n^{3/2}t\rfloor})\right)_{t\geq 0}\rightarrow \phi({X}^w)$ in distribution in the space $D(\mathbb{R}_+,l^1)$, conditional on $X^n_0=\rho_n$ and $X^w_0=\rho_w$.
\end{itemize}}

Note that, in \cite{Croydoncbp}, in place of $\nu^n$, the uniform measure on the vertices of $\mathcal{T}_n$, $\mu^n$ say, was considered. However, it is easy to check that the Prohorov distance between $(2n)^{-1}\nu^n(\phi_n^{-1}(n^{1/2}\cdot))$ and $n^{-1}\mu^n(\phi_n^{-1}(n^{1/2}\cdot))$ is bounded above by $n^{-1/2}$, and so \cite{Croydoncbp}, Theorem 1.1 does indeed imply the above convergence result for measures.

Whenever $n^{-1/2}w_n\rightarrow w\in\mathcal{W}^*$, it is a simple consequence of Proposition \ref{prop} to check that if $X^w$ admits a transition density $(p_t(x,y))_{x,y\in\mathcal{T}_w,t>0}$ which is jointly continuous in $(t,x,y)$, then Assumption 1 holds with $E=l^1$, $F=\phi(\mathcal{T}_w)$, $\rho=0$, $\nu=\nu_w\circ\phi^{-1}$, $q_t(x)=p_t(\rho,x)$, $G^n=n^{-1/2}\phi_n(\mathcal{T}_n)$ and $\alpha(n)=n^{1/2}$, $\beta(n)=2n$, $\gamma(n)=n^{3/2}$. However, by applying \cite{Croydoncrt}, Theorem 6.2 and Corollary 8.5, it is possible to assume that $X^w$ does indeed admit a suitable transition density for $w\in \mathcal{W}^*$. Furthermore, since $\mathcal{T}_n$ is a graph tree for each $n$, it is a fact that the resistance metric $R_{\mathcal{T}_n}$ is identical to the usual graph distance $d_{\mathcal{T}_n}$, and therefore Assumption 5 holds with $\kappa=1$. Hence the following result is true, where we use the notation $q^n:=q^{\mathcal{T}_n}$.

{\thm\label{treeq} There exists a set $\mathcal{W}^*\subseteq C([0,1],\mathbb{R}_+)$ such that $W\in \mathcal{W}^*$ almost-surely, and if $n^{-1/2}w_n\rightarrow w$ in $C([0,1],\mathbb{R}_+)$ for some $w\in\mathcal{W}^*$, then there exists, for each $n$, an isometric embedding $\phi_n$ of $(V(\mathcal{T}_n),d_{\mathcal{T}_n})$ into $l^1$ and also an isometric embedding $\phi$ of $(\mathcal{T}_w,d_{\mathcal{T}_w})$ into $l^1$ such that: in addition to the convergence results of Proposition \ref{prop}, for every compact interval $I\subset (0,\infty)$,
\[\lim_{n\rightarrow\infty}\sup_{x\in\mathcal{T}_w}\sup_{t\in I}\left|2nq^{n}_{\lfloor n^{3/2}t\rfloor}(\tilde{g}_n(x))-p_t(\rho_w,x)\right|=0,\]
where, for $x\in\mathcal{T}_w$, $\tilde{g}_n(x)$ is a point in $V(\mathcal{T}_n)$ minimising the $l^1$-distance between $\phi(x)$ and $n^{-1/2}\phi_n(y)$ over $y\in V(\mathcal{T}_n)$.}
\bigskip

We now present a topology for transition densities on graphs and metric spaces that allows us to state a version of this result that does not involve the underlying metric space $E=l^1$. A particular motivation for doing this is that it allows us to deduce, in addition to the above quenched local limit theorem, a corresponding distributional result.

For an interval $I\subseteq[0,\infty)$, let $\tilde{\mathcal{M}}_I$ be the collection of triples of the form $(F,\rho_F,q^F)$, where $F=(F,d_F)$ is a non-empty compact metric space, $\rho_F$ is a distinguished element of $F$ and $q^F=(q^F_t(x))_{x\in F, t>0}$ is a jointly continuous real-valued function of $(t,x)$. We say two elements, $(F,\rho_F,q^F)$ and $(F',\rho_{F'},q^{F'})$, of $\tilde{\mathcal{M}}_I$ are equivalent if there exists an isometry $f:F\rightarrow {F'}$ such that $f(\rho_F)=\rho_{F'}$ and $q_t^{F'}\circ f=q_t^F$ for every $t\in I$. Define $\mathcal{M}_I$ to be the set of equivalence classes of $\tilde{\mathcal{M}}_I$ under this relation. Note that we will abuse notation and identify an equivalence class in ${\mathcal{M}}_I$ with a particular element of it. Similarly to the distance between pairs of ``spatial trees'' defined in \cite{LegallDuquesne}, we introduce a distance on $\mathcal{M}_I$ that uses the notion of a correspondence between metric spaces, where, if $F$ and ${F'}$ are two compact metric spaces, a correspondence between $F$ and ${F'}$ is a subset $\mathcal{C}$ of $F\times {F'}$ such that for every $x\in F$ there exists at least one $y\in {F'}$ such that $(x,y)\in\mathcal{C}$ and conversely for every $y\in {F'}$ there exists at least one $x\in F$ such that $(x,y)\in \mathcal{C}$. The distortion of the correspondence $\mathcal{C}$ is defined by $\mathrm{dis}(\mathcal{C}):=\sup\{|d_F(x_1,x_2)-d_{F'}(y_1,y_2)|:(x_1,y_1),(x_2,y_2)\in\mathcal{C}\}$, and we set, for $(F,\rho_F,q^F),({F'},\rho_{F'},q^{F'})\in\mathcal{M}_I$,
\[\Delta_I\left((F,\rho_F,q^F),({F'},\rho_{F'},q^{F'})\right):=\inf_{\substack{\mathcal{C}\in\mathfrak{C}(F,{F'}):\\ (\rho_F,\rho_{F'})\in\mathcal{C}}}\delta_I(\mathcal{C}),\]
where $\delta_I(\mathcal{C}):=(\mathrm{dis}(\mathcal{C})+\sup_{(x,y)\in\mathcal{C}, t\in I}\left|q^F_t(x)-q^{F'}_t(y)\right|)$, and $\mathfrak{C}(F,{F'})$ is the set of all correspondences between $F$ and $F'$.

{\lem For any compact interval $I\subseteq [0,\infty)$, $(\mathcal{M}_I,\Delta_I)$ is a separable metric space.}
\begin{proof} That $\Delta_I$ is a metric can be demonstrated by applying a straightforward adaptation of the proof of \cite{BBI}, Theorem 7.3.30, which demonstrates the analogous result for the Gromov-Hausdorff distance between compact metric spaces. To prove separability, first let $\mathcal{F}$ be the countable collection of metric spaces $(F,d_F)$ such that $F$ is a finite set and $d_F$ takes values in $\mathbb{Q}$, and choose a sequence $(i_n)_{n\geq 1}$ that is dense in $I$. For each $n\geq1$, define $\mathcal{M}_I^n$ to be set of equivalence classes of triples of the form $(F,\rho_F,q^F)$ such that $F\in \mathcal{F}$ and, for every $x\in F$, $q_t^F(x)$ takes values in $\mathbb{Q}$ for $t\in \{i_m:m\leq n\}$ and is a linear function of $t$ between the values in $\{i_m:m\leq n\}$. It is an elementary exercise to show that $\cup_{n\geq 1} \mathcal{M}_I^n$ is dense in $(\mathcal{M}_I,\Delta_I)$. Hence, because $\cup_{n\geq 1} \mathcal{M}_I^n$ is countable, $(\mathcal{M}_I,\Delta_I)$ is separable.
\qed\end{proof}

The following result contains a distributional analogue of Theorem \ref{treeq} that applies the above topology. Note that for a graph $G$, we extend the discrete time function $(q^G_m(x))_{x\in V(G),m\geq 0}$ to continuous time by linear interpolation at each vertex. Thus we can view $((V(G),d_G),\rho(G),(q^G_t(x))_{x\in V(G),t\in I})$ as an element of $\mathcal{M}_I$ for every finite graph $G$.

{\thm Fix a compact interval $I\subset(0,\infty)$. Suppose that $(\mathcal{T}_n)_{n\geq2}$ is a sequence of random rooted ordered graph trees whose search-depth functions $(w_n)_{n\geq 2}$ converge in distribution to $W$, the Brownian excursion normalised to have length one, then
\[\left(\left(V(\mathcal{T}_n),n^{-1/2}d_{\mathcal{T}_n}\right),\rho_n,\left(2nq^n_{n^{3/2}t}(x)\right)_{x\in V(\mathcal{T}_n),t\in I}\right)\]
converges in distribution to \[\left(\left(\mathcal{T}_W,d_{\mathcal{T}_W}\right),\rho_W,\left(p_{t}(\rho_W,x)\right)_{x\in \mathcal{T}_W,t\in I}\right)\]
in the space $(\mathcal{M}_I,\Delta_I)$.}
\begin{proof} By the separability of $C([0,1],\mathbb{R}_+)$, it is possible to assume that we have realisations of $(\mathcal{T}_n)_{n\geq2}$ and $W$ such that $n^{-1/2}w_n\rightarrow W$ almost-surely. It is an easy consequence of Theorem \ref{treeq} (and Assumption 2) that
\[\left(\left(V(\mathcal{T}_n),n^{-1/2}d_{\mathcal{T}_n}\right),\rho_n,\left(2nq^n_{n^{3/2}t}(x)\right)_{x\in V(\mathcal{T}_n),t\in I}\right)\]
converges to \[\left(\left(\mathcal{T}_W,d_{\mathcal{T}_W}\right),\rho_W,\left(p_{t}(\rho_W,x)\right)_{x\in \mathcal{T}_W,t\in I}\right)\]
in the space $(\mathcal{M}_I,\Delta_I)$ almost-surely, and the result follows.
\qed\end{proof}

\subsection{Local homogenisation for nested fractals}\label{nf}

We start this section by introducing unbounded nested fractal sets, which will later appear as scaling limits of the associated fractal graphs. Suppose $(\psi_i)_{i=1}^N$ is a family of $L^{-1}$-similitudes on $\mathbb{R}^d$ for some $L>1$, by which we mean that, for each $i$, $\psi_i$ is a map from $\mathbb{R}^d$ to $\mathbb{R}^d$ that satisfies $|\psi_i(x)-\psi_i(y)|=L^{-1}|x-y|$, for every $x,y\in \mathbb{R}^d$, where $|\cdot-\cdot|$ is the usual Euclidean distance on $\mathbb{R}^d$. We assume that the collection $(\psi_i)_{i=1}^N$ satisfies the open set condition; this means that there exists a non-empty bounded set $O\subseteq\mathbb{R}^d$ such that $(\psi_i(O))_{i=1}^N$ are disjoint and $\cup_{i=1}^N\psi_i(O)\subseteq O$. Since $(\psi_i)_{i=1}^N$ is a family of contraction maps, there exists a unique non-empty compact set $K$ such that $K=\cup_{i=1}^N \psi_i(K)$, which we will suppose is connected. Write the set of fixed points of $(\psi_i)_{i=1}^N$ as $\Xi$, and define the collection of essential fixed points of $(\psi_i)_{i=1}^N$ by \[V_0:=\left\{x\in \Xi:\mbox{ $\exists i,j\in\{1,\dots,N\}$, $i\neq j$ and $y\in \Xi$ such that $\psi_i(x)=\psi_i(y)$}\right\}.\]
Throughout, we assume that $\# V_0\geq 2$. The compact set $K$ is then said to be a nested fractal if it satisfies the following finite ramification and symmetry properties.
\begin{itemize}
  \item If $i_1\dots i_n$ and $j_1\dots j_n$ are distinct sequences in $\{1,\dots,N\}$, then
  \[\psi_{i_1\dots i_n}(K)\cap\psi_{j_1\dots j_n}(K)=\psi_{i_1\dots i_n}(V_0)\cap\psi_{j_1\dots j_n}(V_0),\]
  where $\psi_{i_1\dots i_n}:=\psi_{i_1}\circ\dots\circ\psi_{i_n}$.
  \item If $x,y\in V_0$, then the reflection in the hyperplane $H_{xy}:=\{z\in \mathbb{R}^d:\:|z-x|=|z-y|\}$ maps $V_n$ to itself, where
  \[V_n:=\bigcup_{i_1,\dots,i_n=1}^N \psi_{i_1\dots i_n}(V_0).\]
\end{itemize}
Without loss of generality, we assume that $\psi_1(x)=L^{-1}x$ and $0\in V_0$. The unbounded nested fractal which will be of interest in this section is then defined by
\[F:=\bigcup_{n\geq 0}L^n K.\]
Although the embedding of $F$ into Euclidean space has been important for its construction, it will not be particularly important in what follows. Instead, we consider an intrinsic geodesic metric $d_F$ on $F$, as defined in Section 3 of \cite{FHK} (we assume that the size vector introduced there is simply $\tilde{\mathbf{r}}=(1,\dots,1)$), which satisfies the properties presented in the following lemma. In particular, we describe how the metric $d_F$ is related to both the Euclidean metric in $\mathbb{R}^d$, and the shortest path graph distance on the graph $\mathcal{G}$, defined by setting
\[V(\mathcal{G}):=\bigcup_{n\geq 0}L^n V_0\]
and edge set $E(\mathcal{G})$ equal to
\[\left\{\{L^n\psi_{i_1\dots i_n}(x),L^n\psi_{i_1\dots i_n}(y)\}:\:x,y\in V_0,\:x\neq y,\:i_1,\dots,i_n\in\{1,\dots,N\},\:n\geq 0\right\}.\]
We also record a scaling formula for $d_F$, which is proved in \cite{FHK}.

{\lem There exists a metric $d_F$ on $F$, which satisfies the midpoint property, and moreover, there exist constants $c_1,c_2,c_3,c_4\in(0,\infty)$ and $\alpha\in(1,\infty)$ such that
\begin{equation}\label{chemical}
c_1d_F(x,y)\leq |x-y|^{d_c}\leq c_2 d_F(x,y),\hspace{20pt}\forall x,y \in F,
\end{equation}
where $d_c:=\ln \alpha/\ln L$, and also
\begin{equation}\label{df}c_3d_F(x,y)\leq d_\mathcal{G}(x,y)\leq c_4d_F(x,y),\hspace{20pt}\forall x,y \in V(\mathcal{G}).\end{equation}
Finally, $d_F$ can be constructed so that
\begin{equation}\label{scale}d_F(L x,L y)=\alpha d_F(x,y),\hspace{20pt}\forall x,y \in F.\end{equation}}

Observe that this result implies all the conditions on $(F,d_F)$ that are required in the introduction, where we suppose throughout the remainder of this section that $(E,d_E):=(F,d_F)$ and $\rho:=0$. Furthermore, it is a standard result that $K$ has Hausdorff dimension $\ln N/\ln L$ with respect to the Euclidean metric (see \cite{Falconerbook}, for example), and the same is true of $F$. We will denote by $\nu$ the $(\ln N/\ln L)$-dimensional Hausdorff measure on $F$ with respect to the Euclidean metric, and note that it is easy to check (by applying (\ref{chemical}) and the symmetries of the fractal) that there exist constants $c_1,c_2\in (0,\infty)$ such that $c_1 r^{d_f}\leq \nu(B_F(x,r))\leq c_2r^{d_f}$, for every $x\in F$ and $r>0$, where $d_f:=\ln N/\ln\alpha$ and $B_F(x,r)$ is the open ball with centre $x$ and radius $r$ in $(F,d_F)$. Hence $\nu$ satisfies the properties required of the measure on the metric space $(F,d_F)$ in the introduction. The following continuity result will also be useful.

{\lem \label{nucont} For every $x\in F$ and $r>0$, $\nu(\partial B_F(x,r))=0$, where $\partial B_F(x,r):=\overline{B}_F(x,r)\backslash B_F(x,r)$.}
\begin{proof} We prove the corresponding result for $K$, the lemma then follows by rescaling. As a straightforward consequence of \cite{FHK}, Proposition 3.6, there exist constants $c_1,c_2\in (0,\infty)$ such that
\begin{equation}\label{fhk}
c_1\alpha^{-n}< \sup_{x\in \psi_{i_1\dots i_n}(K)} \inf_{y\in K\backslash\psi_{i_1\dots i_n}(K)} d_K(x,y)\leq \sup_{x,y\in \psi_{i_1\dots i_n}(K)} d_K(x,y)< c_2\alpha^{-n}
\end{equation}
for every $i_1,\dots,i_n\in\{1,\dots,N\}$, $n\in \mathbb{N}$, where $d_K:=d_F|_{K\times K}$. Choose $M$ to be an integer strictly greater than $\ln(8c_2/c_1)/\ln\alpha$.

Let $x\in F$, $r>0$ and let $n_0$ be an integer chosen to satisfy $r>2c_2\alpha^{-n_0M}$. We now claim that if $\mathcal{I}_n\subseteq \{1,\dots,N\}^{nM}$ is chosen so that $(\psi_{i_1\dots i_{nM}}(K))_{i_1\dots i_{nM}\in \mathcal{I}_n}$ is a cover for $\partial B_K(x,r)$ and $n\geq n_0$, then there exists a set $\mathcal{I}_{n+1}\subseteq \{1,\dots,N\}^{(n+1)M}$ for which $(\psi_{i_1\dots i_{(n+1)M}}(K))_{i_1\dots i_{(n+1)M}\in \mathcal{I}_{n+1}}$ is a cover for $\partial B_K(x,r)$ and $\#\mathcal{I}_{n+1}\leq (N^M-1)\#\mathcal{I}_n$. Let $(i_1,\dots,i_{nM})\in \mathcal{I}_n$. Clearly, we can assume that there exists an $x_0\in \psi_{i_1\dots i_{nM}}(K)$ such that $d_K(x,x_0)=r$ (if not, then we can discard $(i_1,\dots,i_{nM})$ from $\mathcal{I}_n$), and, by (\ref{fhk}), there exists an $x_1\in \psi_{i_1\dots i_{nM}}(K)$ such that $B_K(x_1,c_1\alpha^{-nM})\subseteq \psi_{i_1\dots i_{nM}}(K)$. We have, applying (\ref{fhk}) and our choice of $n_0$,
\[d_K(x,x_1)\geq |d_K(x,x_0)-d_K(x_0,x_1)|\geq r- c_2\alpha^{-nM}\geq c_2\alpha^{-nM}.\]
Thus there exists an $x_2\in\psi_{i_1\dots i_{nM}}(K)$ on the geodesic path from $x$ to $x_1$ that satisfies $d_K(x,x_2)=d_K(x,x_1)-\tfrac{1}{2}c_1\alpha^{-nM}$. It immediately follows that we can find $j_1,\dots,j_M\in\{1,\dots,N\}$ and $x_3\in\psi_{i_1\dots i_{nM}j_1\dots j_M}(K)$ such that $|r-d_K(x,x_3)|>\tfrac{1}{8}c_1\alpha^{-nM}$. If $z\in \psi_{i_1\dots i_{nM}j_1\dots j_M}(K)$, then
\[|r-d_K(x,z)|\geq |r-d_K(x,x_3)|-d_K(x_3,z)> \tfrac{1}{8}c_1\alpha^{-nM}-c_2\alpha^{-(n+1)M}>0.\]
In particular, this implies that $\psi_{i_1\dots i_{nM}j_1\dots j_M}(K)\cap \partial B_K(x,r)=\emptyset$. The claim can easily be obtained from this result.

By applying the conclusion of the previous paragraph and the scaling relation at (\ref{chemical}), an elementary argument can be applied to deduce that the Hausdorff dimension of $\partial B_K(x,r)$ with respect to the Euclidean metric is no greater than $\ln (N^M-1)/\ln L^M$. Hence, since $\nu$ was defined as the $\ln N/\ln L$-dimensional Hausdorff measure on $K$, we must have that $\nu(\partial B_K(x,r))=0$ as desired.
\qed\end{proof}

We will henceforth suppose that the weights on the edges in  $E(\mathcal{G})$ are selected randomly from a law $\mathbb{P}_\mu$ on $(0,\infty)^{E(\mathcal{G})}$ which satisfies uniform boundedness and cell independence. By uniform boundedness, we mean that there exist deterministic constants $c_1,c_2\in (0,\infty)$ such that, $\mathbb{P}_\mu$-a.s.,
\begin{equation}\label{ub}
c_1\leq \mu^{\mathcal{G}}_{xy}\leq c_2,\hspace{20pt}\forall \{x,y\}\in E(\mathcal{G}),
\end{equation}
and define cell independence to be the property that for each $n\geq 0$, the collections
\[\left\{\left(\mu_{(L^n\psi_{i_1\dots i_n}(x))(L^n\psi_{i_1\dots i_n}(y))}^\mathcal{G}\right)_{x,y\in V_0,x\neq y}\right\}_{i_1,\dots,i_n\in \{1,\dots,N\}}\]
are independent and have the same distribution as $(\mu_{xy}^\mathcal{G})_{x,y\in V_0,x\neq y}$. Note that we still require $\mu_{xy}^\mathcal{G}=\mu_{yx}^\mathcal{G}$ for every $x,y\in V(\mathcal{G})$, and $\mu_{xy}^\mathcal{G}=0$ if $\{x,y\}\not\in E(\mathcal{G})$. In the next lemma we deduce that measure $\nu^{\mathcal{G}}$ on $V(\mathcal{G})$ associated with such a family of random weights can be rescaled to obtain the  measure $\nu$ on $F$.

{\lem \label{vague} If we denote $\nu^n:=\nu^\mathcal{G}(L^n\cdot)$, then there exists a deterministic constant $c\in (0,\infty)$ such that, $\mathbb{P}_{\mu}$-a.s., the measures $cN^{-n}\nu^n$ converge to $\nu$ in the vague topology on locally finite Borel measures on $(F,d_F)$.}
\begin{proof} By definition, we have that
\begin{equation}\label{decomp}
\nu^n(K):=\sum_{\substack{x,y\in V(\mathcal{G})\cap L^nK:\\\{x,y\}\in E(\mathcal{G})}}\mu_{xy}^\mathcal{G}+\sum_{\substack{x\in V(\mathcal{G})\cap L^nK, \:y\in V(\mathcal{G})\backslash L^nK:\\\{x,y\}\in E(\mathcal{G})}}\mu_{xy}^\mathcal{G}.
\end{equation}
Applying the independence properties of $(\mu_{xy}^\mathcal{G})_{x,y\in V_0,x\neq y}$, the first term is equal in distribution to $\sum_{i_1,\dots,i_n=1}^{N}\xi_{i_1\dots i_n}$, where $(\xi_{i_1\dots i_n})_{i_1,\dots,i_n=1}^N$ are independent copies of $\xi:=\sum_{x,y\in V_0}\mu_{xy}^\mathcal{G}$. It follows that, for every $\varepsilon>0$,
\[\mathbb{P}_\mu\left(\left|N^{-n}\sum_{i_1,\dots,i_n=1}^{N}\xi_{i_1\dots i_n}-\mathbb{E}_\mu\xi\right|>\varepsilon\right)\leq \frac{\mathbb{E}_\mu\xi^2}{N^n\varepsilon^2}.\]
Thus, by applying the Borel-Cantelli lemma, we are able to deduce that, when multiplied by $N^{-n}$, the first term of (\ref{decomp}) converges to $\mathbb{E}_\mu\xi$, $\mathbb{P}_\mu$-a.s.

Under the assumption of uniform boundedness, the second term of (\ref{decomp}) is bounded above deterministically by $c_1\#\{\{x,y\}\in E(\mathcal{G}):\:x\in V(\mathcal{G})\cap L^nK,\:y\in V(\mathcal{G})\backslash L^nK\}$ for some constant $c_1$. It is straightforward to check that $\mathcal{G}$ is a graph of bounded degree (cf. \cite{Barlow}, Proposition 5.21) and by combining this fact with the finite ramification property of $K$, it is possible to deduce that $\#\{\{x,y\}\in E(\mathcal{G}):\:x\in V(\mathcal{G})\cap L^nK,\:y\in V(\mathcal{G})\backslash L^nK\}$ is bounded by a constant that is independent of $n$. This completes the proof that $c_2N^{-n}\nu^n(K)\rightarrow \nu(K)$, $\mathbb{P}_\mu$-a.s., where $c_2:=\nu(K)/\mathbb{E}_\mu\xi$.

Applying the self-similarity of $F$, the above argument is easily generalised to yield, $\mathbb{P}_\mu$-a.s.,
\[\lim_{n\rightarrow\infty}c_2N^{-n}\nu^n\left(L^{m_1}\psi_{i_1\dots i_{m_2}}(K)\right)= \nu \left(L^{m_1}\psi_{i_1\dots i_{m_2}}(K)\right),\]
for every $i_1,\dots,i_{m_2}\in\{1,\dots,N\}$, $m_1,m_2\in\mathbb{N}$. The lemma follows from this by applying \cite{Bill2}, Theorem 2.3, for example.
\qed\end{proof}

Our description of the transition density asymptotics and scaling limit of the simple random walk on $\mathcal{G}$ will be presented in terms a resistance-scaling factor $\lambda$, which appears as an ``eigenvalue'' for the renormalisation map that we now introduce. For a set of non-negative weights $(C_{xy})_{x,y\in V^0, x\neq y}$ which satisfy $C_{xy}=C_{yx}$, a quadratic form $\mathcal{E}_C$ on $\mathbb{R}^{V_0}$ can be constructed by setting
\[\mathcal{E}_C(f,f):=\sum_{x,y\in V_0,x\neq y}C_{xy}(f(x)-f(y))^2.\]
Replicating this form $N$ times, we set
\[\mathcal{E}^1_C(f,f):=\sum_{i=1}^N \mathcal{E}_C(f\circ\psi_i,f\circ \psi_i),\]
which defines a quadratic form on $\mathbb{R}^{V_1}$. Now, restrict this form to $V_0$ using the trace operator, as defined by $\mathrm{Tr}(\mathcal{E}^1_C|V^0)(f,f):=\inf\{\mathcal{E}_C^1(g,g):\:g|_{V_0}=f\}$. The resulting operator $\mathrm{Tr}(\mathcal{E}^1_C|V^0)$ is of the form $\mathcal{E}_{\Lambda(C)}$ for non-negative weights $(\Lambda(C)_{xy})_{x,y\in V^0, x\neq y}$ which satisfy $\Lambda(C)_{xy}=\Lambda(C)_{yx}$. It is known that there exists a non-degenerate fixed point to the map $C\mapsto \Lambda(C)$ which satisfies $\Lambda(C)=\lambda^{-1}C$, for some $\lambda>0$ (see \cite{Barlow}, Theorem 6.23 for example). In fact, $\lambda$ is uniquely determined, so that it is the same for any non-degenerate fixed point (\cite{Barlow}, Corollary 6.20). Moreover, we can also assume that $\lambda>1$ (\cite{Barlow}, Corollary 6.28).

In the subsequent lemma, we summarise results for the continuous time simple random walk $Y^\mathcal{G}$ on the graph $\mathcal{G}$, its transition density $(\tilde{p}^\mathcal{G}_t(x,y))_{x,y\in V(\mathcal{G}),t>0}$ and the resistance metric $R_\mathcal{G}$  determined from the corresponding Dirichlet form by the formula at (\ref{resdef}). The constant $d_w$ is defined to be equal to $\ln (N\lambda)/\ln\alpha$, and $d_f:=\ln N/\ln \alpha$, as above.

{\lem \label{nfsrw} $\mathbb{P}_\mu$-a.s., there exist (random) constants $c_1,c_2,c_3,c_4,c_5,c_6\in(0,\infty)$ such that
\begin{equation}\label{nfghk}
\tilde{p}^\mathcal{G}_t(x,y)\leq c_1 t^{-d_f/d_w} \exp\left(-c_2\left(\frac{d_\mathcal{G}(x,y)^{d_w}}{t}\right)^{1/(d_w-1)}\right),\hspace{20pt}\forall x,y\in V(\mathcal{G}),t>0,
\end{equation}
\begin{equation}\label{exit}\tilde{\mathbf{P}}_x^\mathcal{G}\left(\tau^\mathcal{G}(x,r)\leq t\right)\leq c_3\exp\left(-c_4\left(\frac{r^{d_w}}{t}\right)^{1/(d_w-1)}\right),\hspace{20pt}\forall x\in V(\mathcal{G}),t,r>0,
\end{equation}
where $\tau^\mathcal{G}(x,r):=\inf\{t>0:\:Y^\mathcal{G}_t\not\in B_{\mathcal{G}}(x,r)\}$ is the exit time of the simple random walk $Y^\mathcal{G}$ from the graph ball $B_{\mathcal{G}}(x,r)$, and also
\begin{equation}\label{rescont}
c_5 d_\mathcal{G}(x,y)^\kappa \leq R_\mathcal{G}(x,y)\leq c_6 d_\mathcal{G}(x,y)^\kappa,\hspace{20pt}\forall x,y\in V(\mathcal{G}),
\end{equation}
where $\kappa:=\ln \lambda/\ln\alpha$.}
\begin{proof} For the discrete time simple random walk on the unweighted nested fractal graph $\mathcal{G}$ (i.e. $\mu_{xy}^\mathcal{G}=1$ for every $\{x,y\}\in E(\mathcal{G})$), the parabolic Harnack inequality with exponent $d_w$ is known to hold (see \cite{GT}, Theorem 3.1 and \cite{HamKum2}, Corollary 4.13). Hence, by \cite{HamKum2}, Theorem 5.11, the same is true for any uniformly bounded set of weights $\mu^\mathcal{G}$, $\mathbb{P}_\mu$-a.s. Given this property, the discrete time versions of (\ref{nfghk}) and (\ref{exit}), as well as (\ref{rescont}), are an application of results appearing in \cite{BCK}. Similar arguments can be used to prove the corresponding continuous time results (alternatively, once (\ref{rescont}) is established, the arguments of \cite{BCK} can be adapted to continuous time directly to yield (\ref{nfghk}) and (\ref{exit})).
\qed\end{proof}

We can use this lemma to prove a $\mathbb{P}_\mu$-a.s. tightness result for the law of the continuous time simple random walk $Y^\mathcal{G}$. Note that, for $x\in F$, we write $\tilde{g}_n(x)$ to represent the point in $V(\mathcal{G})$ closest to $L^n x$.

{\lem \label{tightness}For every compact interval $I\subset (0,\infty)$, $x,y\in F$ and $r>0$, we have that, $\mathbb{P}_\mu$-a.s.,
\begin{eqnarray*}\lefteqn{\lim_{\delta\rightarrow 0}\limsup_{n\rightarrow \infty}\sup_{\substack{s,t\in I\\|s-t|< \delta}}\hspace{-.2pt}\left|\tilde{\mathbf{P}}^\mathcal{G}_{\tilde{g}_n(x)}\left(L^{-n}Y^\mathcal{G}_{(N\lambda)^n s}\in B_F(y,r)\right)\right.}\\
&&\hspace{100pt}\left.-\tilde{\mathbf{P}}^\mathcal{G}_{\tilde{g}_n(x)}\left(L^{-n}Y^\mathcal{G}_{(N\lambda)^n t}\in B_F(y,r)\right)\right|=0.
\end{eqnarray*}}
\begin{proof} Fix a compact interval $I\subset (0,\infty)$, $x,y\in F$ and $r,\varepsilon>0$, and write $B=B_F(y,r)$. Some elementary analysis allows us to conclude that, for any $\eta>0$,
\begin{eqnarray}
\lefteqn{\sup_{\substack{s,t\in I\\|s-t|< \delta}}\left|\tilde{\mathbf{P}}^\mathcal{G}_{\tilde{g}_n(x)}\left(Z_s^n\in B\right)-\tilde{\mathbf{P}}^\mathcal{G}_{\tilde{g}_n(x)}\left(Z_t^n\in B\right)\right|}\nonumber\\
&\leq & 2\sup_{\substack{s,t\in I\\0<t-s<\delta}}\tilde{\mathbf{P}}^\mathcal{G}_{\tilde{g}_n(x)}\left(d_F(Z_s^n,Z_t^n)>\eta\right)+
2\sup_{t\in I}\tilde{\mathbf{P}}^\mathcal{G}_{\tilde{g}_n(x)}\left(Z_t^n\in B_\eta\backslash B\right)\label{twoterms},
\end{eqnarray}
where $B_\eta:=B_F(y,r+\eta)$ and we denote $L^{-n}Y^\mathcal{G}_{(N\lambda)^n t}$ by $Z^n_t$ in this proof. For the second term, we can apply the heat kernel bound of (\ref{nfghk}) and the measure convergence of Lemma \ref{vague} to deduce that, $\mathbb{P}_\mu$-a.s., there exists a finite constant $c_1$ such that, for every $\eta>0$
\[\limsup_{n\rightarrow \infty}2\sup_{t\in I}\tilde{\mathbf{P}}^\mathcal{G}_{\tilde{g}_n(x)}\left(Z_t^n\in B_\eta\backslash B\right)\leq c_1\nu(B_\eta\backslash B).\]
By Lemma \ref{nucont}, this upper bound is less than $\varepsilon$ for suitably small $\eta$.

For the first term in (\ref{twoterms}), we apply the Markov property of $Y^\mathcal{G}$ and the metric approximation result of (\ref{df}) to obtain the existence of a (deterministic) non-zero constant $c_2$ such that
\begin{eqnarray*}
2\sup_{\substack{s,t\in I\\0<t-s<\delta}}\tilde{\mathbf{P}}^\mathcal{G}_{\tilde{g}_n(x)}\left(d_F(Z_s^n,Z_t^n)>\eta\right)&\leq & 2\sup_{t\in[0,\delta)}\sup_{z\in V(\mathcal{G})} \tilde{\mathbf{P}}^\mathcal{G}_z\left(d_\mathcal{G}(z,Y^\mathcal{G}_{(N\lambda)^nt})>c_2\alpha^n\eta\right)\\
&\leq & 2\sup_{z\in V(\mathcal{G})} \tilde{\mathbf{P}}^\mathcal{G}_z\left(\tau^\mathcal{G}(z,c_2\alpha^n \eta)\leq (N\lambda)^n\delta\right),
\end{eqnarray*}
where $\tau^\mathcal{G}(\cdot,\cdot)$ is the exit time defined in Lemma \ref{nfsrw}. Consequently, the upper bound for the exit time distribution at (\ref{exit}) implies that, $\mathbb{P}_\mu$-a.s.,
\[\lim_{\delta\rightarrow 0}\limsup_{n\rightarrow\infty}2\sup_{\substack{s,t\in I\\0<t-s<\delta}}\tilde{\mathbf{P}}^\mathcal{G}_{\tilde{g}_n(x)}\left(d_F(Z_s^n,Z_t^n)>\eta\right)=0.\]
In combination with the conclusion of the previous paragraph, this completes the proof.
\qed\end{proof}

We continue by describing how the homogenisation result of \cite{homog} can be applied in our situation. Set, for $f\in \mathbb{R}^{V_n}$,
\begin{eqnarray}
\lefteqn{\mathcal{E}_\mu^n(f,f):=}\label{nth}\\
&&\sum_{i_1,\dots,i_n=1}^N\sum_{x,y\in V_0, x\neq y}\mu^\mathcal{G}_{(L^{n}\psi_{i_1\dots i_n}(x))(L^{n}\psi_{i_1\dots i_n}(y))}\left(f(\psi_{i_1\dots i_n}(x))-f(\psi_{i_1\dots i_n}(y))\right)^2.\nonumber
\end{eqnarray}
From this, we define $\Lambda^n(\mu)=(\Lambda^n(\mu)_{xy})_{x,y\in V_0,x\neq y}$ to satisfy $\mathcal{E}_{\Lambda^n(\mu)}=\mathrm{Tr}(\mathcal{E}^n_\mu|V_0)$. It is proved in \cite{homog}, Theorem 3.4, that there exists a deterministic $C^\mu=(C^\mu_{xy})_{x,y\in V_0,x\neq y}$ such that
\[\lim_{n\rightarrow\infty}\lambda^n\Lambda^n(\mu)=C^\mu,\]
where the limit is an $L^1$-limit in the space of non-negative weights on the complete graph with vertex set $V_0$. Moreover, $C^\mu$ satisfies $C^\mu_{xy}>0$ for every $x,y\in V_0$, $x\neq y$, and also $\Lambda(C^\mu)=\lambda^{-1} C^\mu$, where $\Lambda$ is the renormalisation map defined above. We will use the weights $C^\mu$ to construct the diffusion on $F$ that arises as the scaling limit of the random walk $Y^\mathcal{G}$ as follows. First, let $\mathcal{E}^n_{C^\mu}$ be a quadratic form on $\mathbb{R}^{V_n}$ which satisfies (\ref{nth}) with $\mu^\mathcal{G}_{(L^{n}\psi_{i_1\dots i_n}(x))(L^{n}\psi_{i_1\dots i_n}(y))}$ replaced by $C^\mu_{xy}$ in each summand, then define
\[\mathcal{E}_K(f,f)=\lim_{n\rightarrow\infty}\lambda^n\mathcal{E}^n_{C^\mu}(f|_{V_n},f|_{V_n})\]
for $f\in \mathcal{F}_K$, where $\mathcal{F}_K:=\{f\in C(K,\mathbb{R}):\:\sup_{n}\lambda^n\mathcal{E}^n_{C^\mu}(f|_{V_n},f|_{V_n})<\infty\}$. It is known (\cite{Barlow}, \cite{Kigami}) that $(\mathcal{E}_K,\mathcal{F}_K)$ is a local, regular (non-degenerate) Dirichlet form on $L^2(K,\nu)$ that satisfies
\[\mathcal{E}_K(f,f)=\lambda\sum_{i=1}^N\mathcal{E}_K(f\circ \psi_i,f\circ\psi_i),\hspace{20pt}\forall f\in \mathcal{F}.\]
For each $n\in \mathbb{Z}$, define a renormalisation operator $\sigma_n$ by setting $\sigma_n(f)(x)=f(L^nx)$ for $x\in K$ and $f:K^n\rightarrow \mathbb{R}$, where $K^n:=L^nK$. If we set $\mathcal{F}_{K^n}:=\sigma_{-n}\mathcal{F}_K$ and
\[\mathcal{E}_{K^n}(f,f):=\mathcal{E}_K(\sigma_n(f),\sigma_n(f)),\hspace{20pt}\forall f\in \mathcal{F}_{K^n},\]
then it is possible to define a local, regular Dirichlet form $(\mathcal{E}_F,\mathcal{F}_F)$ on $L^2(F,\nu)$ by setting
\[\mathcal{E}_F(f,f):=\lim_{n\rightarrow\infty}\lambda^n\mathcal{E}_{K^n}(f|_{K^n},f|_{K^n}),\hspace{20pt}\forall f\in \mathcal{F}_F,\]
where $\mathcal{F}_F$ is the collection of functions $f\in L^2(F,\nu)$ that satisfy $f|_{K^n}\in\mathcal{F}_{K^n}$ for every $n\geq 0$ and $\lim_{n\rightarrow\infty}\mathcal{E}_{K^n}(f|_{K^n},f|_{K^n})<\infty$, see \cite{FHK}, Theorem 2.7. Finally, the associated $\nu$-symmetric diffusion $X=((X_t)_{t\geq 0},\mathbf{P}_x,x\in F)$ admits a transition density $(p_t(x,y))_{x,y\in F,t>0}$ that is jointly continuous in $(t,x,y)$, see \cite{FHK}, Lemma 4.6, and satisfies
\begin{equation}\label{nfhk}
p_t(x,y)\leq c_1 t^{-d_f/d_w}\exp\left(-c_2\left(\frac{d_F(x,y)^{d_w}}{t}\right)^{1/(d_w-1)}\right),\hspace{20pt}\forall x,y\in F,t>0,
\end{equation}
where $d_f=\ln N/\ln \alpha$ and $d_w=\ln(N\lambda)/\ln\alpha$, as before, and $c_1,c_2$ are constants taking values in $(0,\infty)$, see \cite{FHK}, Theorem 5.7. Finally, we have the following important scaling result for the simple random walk $Y^\mathcal{G}$.

{\lem[\cite{homog}, Theorem 7.3] \label{dconv} There exists a deterministic constant $c\in(0,\infty)$ such that, for every $\varepsilon>0$, $x\in F$ and bounded $f\in C(D(\mathbb{R}_+,F),\mathbb{R})$, we have
\[\lim_{n\rightarrow\infty}\mathbb{P}_\mu\left(\left|\tilde{\mathbf{E}}^\mathcal{G}_{\tilde{g}_n(x)}(f(L^{-n}Y^\mathcal{G}_{c(N\lambda)^n\cdot}))-\mathbf{E}_x(f(X_\cdot))\right|>\varepsilon\right)=0,\]
where $\tilde{\mathbf{E}}^\mathcal{G}_{\tilde{g}_n(x)}$ is the expectation under the continuous time simple random walk law $\tilde{\mathbf{P}}_{\tilde{g}_n(x)}^\mathcal{G}$, and $\mathbf{E}_x$ is the expectation under the law $\mathbf{P}_x$ of the Markov process $X$.}
\bigskip

Given all the above information, it is straightforward to deduce properties of the graph $\mathcal{G}$ from which a local limit theorem for nested fractal graphs with random weights can be deduced. We set $G^n:=L^{-n}\mathcal{G}$, by which we mean that $V(G^n)=L^{-n}V(\mathcal{G})$, $E(G^n):=\{\{L^{-n}x,L^{-n}y\}:\:\{x,y\}\in E(\mathcal{G})\}$ and $\mu_{xy}^{G^n}=\mu^\mathcal{G}_{(L^nx)(L^ny)}$.

{\propn\label{nfp} The graphs $(G^n)_{n\geq 0}$ satisfy the following assumptions for $\alpha(n)=\alpha^n$, $\beta(n)=c_1N^n$ and $\gamma(n)=c_2(N\lambda)^n$, for some deterministic constants $c_1,c_2\in(0,\infty)$.\\
(i) Assumption 1(a), 1(b) hold. Assumption 1(c) holds $\mathbb{P}_\mu$-a.s. Furthermore, for every compact interval $I\subset (0,\infty)$, $x,y\in F$,
\begin{equation}\label{label}
\lim_{n\rightarrow\infty}\mathbb{P}_{\mu}\left(\sup_{t\in I}\left|\tilde{\mathbf{P}}_{{g}_n(x)}^{G^n}\left(Y^{G^n}_{\gamma(n)t}\in B_E(y,r)\right)-\mathbf{P}_x\left(X_t\in B_E(y,r)\right)\right|>\varepsilon\right)=0.
\end{equation}
(ii) The continuous time version of Assumption \^{2} holds $\mathbb{P}_\mu$-a.s.\\
(iii) The transition density of $X$ satisfies
\[\lim_{t\rightarrow\infty}\sup_{x,y\in F}{p}_{t}(x,y)=0,\hspace{20pt}\lim_{r\rightarrow\infty}\sup_{x\in B_{F}(\rho,R)}\sup_{y\in F\backslash B_F(\rho,r)} \sup_{t\in I}{p}_{t}(x,y)=0,\]
for any compact interval $I\subset (0,\infty)$, $R>0$. Moreover, $\mathbb{P}_\mu$-a.s.,
\[\lim_{t\rightarrow\infty}\limsup_{n\rightarrow\infty}\sup_{x,y\in V(G^n)} \beta(n)\tilde{p}^{G^n}_{\gamma(n)t}(x,y)=0,\]
and, for any compact interval $I\subset (0,\infty)$, $R>0$,
\[\lim_{r\rightarrow\infty}\limsup_{n\rightarrow\infty}\sup_{x\in B_{G^n}(\rho,\alpha(n)R)}\sup_{y\in V(G^n)\backslash B_{G^n}(\rho,\alpha(n)r)} \sup_{t\in I}\beta(n)\tilde{p}^{G^n}_{\gamma(n)t}(x,y)=0.\]}
\begin{proof} Assumption \ref{first}(a) is a simple consequence of (\ref{df}) and the scaling relation at (\ref{scale}). By construction, $V(G^n)\subseteq V(G^{n+1})$ for every $n\geq 0$ and also $\cup_{n\geq0}V(G^n)$ is dense in $(F,d_F)$, thus Assumption \ref{first}(b) holds. The measure convergence of (\ref{randommeasconv}) is implied by Lemmas \ref{nucont} and \ref{vague}. For the remaining claim of (i), we apply Lemma \ref{dconv} to deduce that there exists a constant $c\in(0,\infty)$ such that, for any $0<t_1<\dots<t_k$, $x,y\in F$, $r,\varepsilon>0$, we have
\[\lim_{n\rightarrow\infty}{\mathbb{P}}_\mu\left(\sup_{i=1,\dots,k}\left| \tilde{\mathbf{P}}^{G^n}_{{g}_n(x)}(Y^{G^n}_{c(N\lambda)^nt_i}\in B_F(y,r))-\mathbf{P}_x(X_{t_i}\in B_F(y,r))\right|>\varepsilon\right)=0.\]
In conjunction with Lemma \ref{tightness}, this implies (\ref{label}).

The control on the resistance metric at (\ref{rescont}) implies that Assumption \ref{fifth} is satisfied $\mathbb{P}_\mu$-a.s. Hence, by the continuous time version of Proposition \ref{52hat}, the continuous time version of Assumption \^{2} holds $\mathbb{P}_\mu$-a.s. To obtain (iii), we apply the heat kernel bounds appearing at (\ref{nfghk}) and (\ref{nfhk}).
\qed\end{proof}

This proposition allows us to obtain the following local limit theorem.

{\thm Fix $T_1, R>0$. Suppose $\mathcal{G}$ is a nested fractal graph with random weights satisfying uniform boundedness and cell independence, then there exist deterministic constants $c_1,c_2\in(0,\infty)$ such that, for every $\varepsilon>0$,
\[\lim_{n\rightarrow\infty}\mathbb{P}_{\mu}\left(\sup_{x\in B_F(\rho,R)}\sup_{y\in F}\sup_{t\geq T_1}\left|c_1N^n\tilde{p}^\mathcal{G}_{c_2(N\lambda)^n t}(\tilde{g}_n(x),\tilde{g}_n(y))-p_t(x,y)\right|>\varepsilon\right)=0.\]}
\begin{proof} By adapting the proof of Theorem \ref{bdd2} to the case of random weights, using similar ideas to those applied in Section \ref{randsec}, and considering the continuous time transition density in place of the discrete time transition density, it is possible to deduce from parts (i) and (ii) of Proposition \ref{nfp} that there exist deterministic constants $c_1,c_2\in(0,\infty)$ such that, for every compact interval $I\subset(0,\infty)$, $R,\varepsilon>0$,
\[\lim_{n\rightarrow\infty}\mathbb{P}_{\mu}\left(\sup_{x,y\in B_F(\rho,R)}\sup_{t\in I}\left|c_1N^n\tilde{p}^\mathcal{G}_{c_2(N\lambda)^n t}(\tilde{g}_n(x),\tilde{g}_n(y))-p_t(x,y)\right|>\varepsilon\right)=0.\]
The theorem easily follows from this by applying the heat kernel decay conditions of Proposition \ref{nfp}(iii).
\qed\end{proof}

Finally, we say that a collection of weights $C=(C_{xy})_{x,y\in V_0,x\neq y}$ is invariant if, for every map $h$ which is a reflection in a hyperplane of the form $H_{xy}$, $x,y\in V_0$, the collection $(C_{h(x)h(y)})_{x,y\in V_0,x\neq y}$ is identical to $(C_{xy})_{x,y\in V_0,x\neq y}$; and it was proved in \cite{Sabot} that there exists a unique non-degenerate invariant set of weights, $C^*$ say, such that $\Lambda(C^*)=\lambda^{-1} C^*$. Thus, if we assume that $(\mu^\mathcal{G}_{xy})_{x,y\in V_0,x\neq y}$ is invariant in distribution (so that $(\mu^\mathcal{G}_{h(x)h(y)})_{x,y\in V_0,x\neq y}$ is equal in distribution to $(\mu^\mathcal{G}_{xy})_{x,y\in V_0,x\neq y}$ for reflections $h$ of the form described), then it follows that $C^\mu=C^*$. The resulting diffusion is known as the Brownian motion on the unbounded nested fractal $F$, and from the above local limit theorem we obtain that if we have a collection of random weights which are invariant in distribution, uniformly bounded and cell independent, then the transition densities of the associated random walk, when rescaled, converge in probability to the transition density of the Brownian motion on the unbounded nested fractal. See Section 7 of \cite{KumKus} for further discussion of invariant weights.

\subsection{Local homogenisation for tree-like Vicsek sets}\label{vicsek}

In this section, we describe a $\mathbb{P}_\mu$-a.s. version of the conclusion of the previous section in a special case. Continuing to apply the notation for nested fractals introduced in Section \ref{nf}, we now assume that $\#V_0=4$ and, moreover, if $\Gamma_n$ is defined to be the graph with vertex set $\{(i_1,\dots,i_n)\}_{i_1,\dots,i_n=1}^N$ and edge set \[E_n:=\left\{(i_1,\dots,i_n),(j_1,\dots,j_n)\}:\:\psi_{i_1\dots i_n}(V_0)\cap\psi_{j_1\dots j_n}(V_0)\neq \emptyset\right\},\]
then $\Gamma_n$ is a graph tree for every $n\in\mathbb{N}$. This class of nested fractals will be referred to as tree-like Vicsek sets, and the Vicsek set (see \cite{HM}, Section 2, for example) is a particular example. The homogenisation problem for tree-like Vicsek sets was studied in \cite{HM}, where the tree-like nature of the graph $\mathcal{G}$ induced by the above assumptions was used to deduce $\mathbb{P}_\mu$-a.s. homogenisation statements, rather than the probabilistic convergence results obtained in \cite{homog} and \cite{KumKus}.

In a slight alteration of our earlier notation, for an edge $e=\{x,y\}\in E(\mathcal{G})$, we now denote $\mu^\mathcal{G}_e:=\mu^\mathcal{G}_{xy}$. The assumptions we will make on the weights are the following: $(\mu_e^\mathcal{G})_{e\in E(\mathcal{G})}$ are independent and identically distributed, have finite second moments and are bounded uniformly below, by which we mean that there exists a constant $c_1>0$ such that $\mu_{e}^\mathcal{G}\geq c_1$ for every $e\in E(\mathcal{G})$, $\mathbb{P}_\mu$-a.s. Under these assumptions, in place of Lemma \ref{dconv}, we have the following. Note that for tree-like Vicsek sets we have $\alpha=\lambda=L$.

{\lem [\cite{HM}, Corollary 1.2] \label{vicsekdconv} There exists a deterministic constant $c\in(0,\infty)$ such that, for every $x\in F$  we have, $\mathbb{P}_\mu$-a.s., if $f\in C(D(\mathbb{R}_+,F),\mathbb{R})$ is bounded, then
\[\lim_{n\rightarrow\infty}\tilde{\mathbf{E}}^\mathcal{G}_{\tilde{g}_n(x)}(f(L^{-n}Y^\mathcal{G}_{c(NL)^n\cdot}))= \mathbf{E}_x(f(X_\cdot)).\]}
\bigskip

Furthermore, we can verify that the continuous time versions of Assumptions \^1 and \^2 hold for tree-like Vicsek sets, where we again consider $G^n=L^{-n}\mathcal{G}$. Note that, unlike the proof of Proposition \ref{nfp}(i), we do not use any transition density estimates.

{\lem $\mathbb{P}_\mu$-a.s., the graphs $(G^n)_{n\geq 0}$ satisfy the continuous time versions of Assumptions \^1 and \^2 for $\alpha(n)=L^n$, $\beta(n)=c_1N^n$ and $\gamma(n)=c_2(NL)^n$, for some deterministic constants $c_1,c_2\in(0,\infty)$.}
\begin{proof} The proof that Assumptions 1(a) and 1(b) hold remains unchanged from Proposition \ref{nfp}. Applying the independence and finite second moments of the weights, it is an elementary exercise to show that the proof of Lemma \ref{vague} can be repeated to deduce Assumption 1(c) in this case. The continuous time version of the convergence at (\ref{2p}) follows from Lemma \ref{vicsekdconv}. Thus Assumption \^1 holds as claimed.

Since the weights are bounded uniformly below, there exists a finite constant $c_1$ such that $R_{G^n}\leq c_1 d_{G^n}$ for every $n$ (see \cite{BCK}, Lemma 2.1); hence the continuous time version of Assumption 5 holds with $\kappa=1$. By the continuous time version of Proposition \ref{52hat}, the continuous time version of Assumption \^2 follows.
\qed\end{proof}

This lemma allows us to apply the continuous time version of Theorem \ref{bdd2} to deduce the subsequent local limit theorem.

{\thm Fix a compact interval $I\subset(0,\infty)$ and $R>0$. Suppose $\mathcal{G}$ is a graph associated with a tree-like Vicsek set equipped with random weights $(\mu_e^\mathcal{G})_{e\in E(\mathcal{G})}$ that are independent and identically distributed, have finite second moments and are bounded uniformly below, then there exist deterministic constants $c_1,c_2\in(0,\infty)$ such that, $\mathbb{P}_\mu$-a.s.,
\[\lim_{n\rightarrow\infty} \sup_{x,y\in B_F(\rho,R)}\sup_{t\in I}\left|c_1N^n\tilde{p}^\mathcal{G}_{c_2(NL)^n t}(\tilde{g}_n(x),\tilde{g}_n(y))-p_t(x,y)\right|=0.\]}

\subsection{Local homogenisation for Sierpinski carpets}\label{carpetsec}

To define generalised Sierpinski carpet graphs, we closely follow \cite{RWSC}. Let $d\geq 2$, $E_0=[0,1]^d$ and $L\in \mathbb{N}$, $L\geq 3$ be fixed. For $n\in\mathbb{Z}$, let $\mathcal{S}_n$ be the collection of closed cubes of side $L^{-n}$ with vertices in $L^{-n}\mathbb{Z}^d$. For $A\subseteq \mathbb{R}^d$, set
\[\mathcal{S}_n(A):=\{S:\:S\subseteq A,\:S\in \mathcal{S}_n\}.\]
For $S\in \mathcal{S}_n$, let $\psi_S$ be the orientation preserving affine map which maps $E_0$ onto $S$. Now suppose that $(\psi_{i})_{i=1}^N$ is a sequence of distinct elements of $(\psi_S)_{S\in\mathcal{S}_{1}}$ and set $E_1=\cup_{i=1}^N\psi_i(E_0)$. We make the following assumptions on $E_1$.
\begin{itemize}
  \item (Symmetry) $E_1$ is preserved by all the isometries of the unit cube $E_0$.  \item (Connectedness) The interior of $E_1$ is connected, and contains a path connecting the hyperplanes $\{x_1=0\}$ and $\{x_1=1\}$.
  \item (Non-diagonality) Let $n\geq 1$ and $B$ be a cube in $E_0$ of side length $2L^{-n}$ which is the union of $2^d$ distinct elements of $\mathcal{S}_n$. Then if the interior of $E_1\cap B$ is non-empty, it is connected.
  \item (Borders included) $E_1$ contains the line segment
  \[\{x:0\leq x_1\leq 1,x_2=\dots=x_d=0\}.\]
\end{itemize}
Note that the non-diagonality assumption stated here can be found in \cite{scu},
and differs from that used in \cite{BarBashd} and \cite{RWSC} for the reasons explained in \cite{scu}\footnote{We thank the anonymous referee for drawing this issue to our attention.}.

Given the maps $(\psi_i)_{i=1}^N$, we can define a generalised Sierpinski carpet $K$ to be the unique non-empty compact set satisfying $K=\cup_{i=1}^N\psi_i(K)$. As in Section \ref{nf}, we denote the associated unbounded carpet $F$. With respect to the Euclidean metric, $F$ has Hausdorff dimension $d_f=\ln N/\ln L$, and we will denote by $\nu$ the $d_f$-dimensional Hausdorff measure on $F$.

To define the corresponding fractal graph, first set
\[P:=\bigcup_{n=1}^\infty L^n\bigcup_{i_1,\dots,i_n=1}^N \psi_{i_1\dots i_n}(E_0),\]
which is the pre-carpet (see \cite{Osada}). Each cube in $\mathcal{S}_0(P)$ has a unique vertex closest to the origin in $\mathbb{R}^d$, let $V(\mathcal{G})$ be the collection of such vertices; in \cite{RWSC}, vertices were chosen to be cube centres, our choice means that $L^{-n}V(\mathcal{G})\subseteq F$ for every $n$. Define $E(\mathcal{G})$ to be the collection of pairs $\{x,y\}$ of elements of $V(\mathcal{G})$ with $|x-y|=1$. The graph of interest in this section will then be $\mathcal{G}=(V(\mathcal{G}),E(\mathcal{G}))$.

Let us now introduce a geodesic metric $d_F$ on $F$. Note that the following result was essentially proved in \cite{Cristea} for the ``standard'' Sierpinski carpet in $\mathbb{R}^2$. As in Section \ref{nf}, for $x\in F$, we write $\tilde{g}_n(x)$ to represent the point in $V(\mathcal{G})$ closest to $L^nx$.

{\lem\label{carpetmetric} For $x,y\in F$, the quantity
\[d_F(x,y):=\lim_{n\rightarrow\infty}L^{-n}d_\mathcal{G}(\tilde{g}_n(x),\tilde{g}_n(y))\]
is well-defined. Moreover $d_F$ is a geodesic metric on $F$ satisfying $d_F(Lx,Lx)=Ld_F(x,y)$ for every $x,y\in F$, and there exists a finite constant $c$ such that
\[|x-y|\leq d_F(x,y)\leq c |x-y|,\hspace{20pt}\forall x,y\in F.\]
Finally, $d_F$ agrees with $d_\mathcal{G}$ on $V(\mathcal{G})$.}
\begin{proof} The proof that $d_F$ is well-defined geodesic metric is similar to \cite{FHK}, Theorem 3.5, and is omitted. The remaining claims are straightforward consequences of the self-similarity and ``borders included'' property of generalised Sierpinski carpets.
\qed\end{proof}

As consequence of this result, if we take $(E,d_E)=(F,d_F)$, $\rho=0$ and $\nu$ as above, then the properties required on the metric spaces and measure in the introduction are satisfied. We can also verify that Assumptions 1(a), 1(b) and 1(c) hold $\mathbb{P}_\mu$-a.s. when, as in Section \ref{nf}, we set $G^n:=L^{-n}\mathcal{G}$ and assume that (using the notation of Section \ref{vicsek}) the weights $(\mu_e^\mathcal{G})_{e\in E(\mathcal{G})}$ are independent, identically-distributed and satisfy uniform boundedness (as at (\ref{ub})).

{\propn \label{carpetprop} $\mathbb{P}_\mu$-a.s., the graphs $(G^n)_{n\geq 0}$ satisfy Assumptions 1(a), 1(b) and 1(c) for $\alpha(n)=L^n$ and $\beta(n)=c_1N^n$, for some deterministic constant $c_1\in(0,\infty)$.}
\begin{proof}Assumptions 1(a) and 1(b) readily follow from the above lemma. To prove Assumption 1(c), we start by demonstrating that $N^n\nu^\mathcal{G}(L^n\tilde{K})$ converges to a deterministic constant $c_2\in(0,\infty)$, $\mathbb{P}_\mu$-a.s., where $\tilde{K}:=\{x\in K:\:x_i\neq 1\mbox{ for any }i=1,\dots,d\}$. Let $e_n$ equal the number of edges of $\mathcal{G}$ that have both end-points in $L^n\tilde{K}$, and ${e}_n'$ represent the number of edges that have exactly one end in $L^n\tilde{K}$, then we have
\[e_{n+1}=Ne_n+e_1\sigma^n,\hspace{20pt}e_n'=d\sigma^n,\]
for every $n\geq 1$, where $\sigma$ is the number of cubes in $\mathcal{S}_1(E_0)$ that lie on a single face of $E_0$. In particular, $\sigma<N$, and so $N^{-n}e_n'\rightarrow 0$ and $N^{-n}e_n\rightarrow c_3$, for some constant $c_3\in (0,\infty)$. Applying the same argument as in the proof of Lemma \ref{vague}, it follows that $N^n\nu^\mathcal{G}(L^n\tilde{K})$ converges as desired. Continuing to imitate the proof of Lemma \ref{vague}, we can extend this to the result that, $\mathbb{P}_\mu$-a.s.,
\begin{equation}\label{vaguecarpet}
c_4N^n\nu^\mathcal{G}(L^n\cdot)\rightarrow\nu
\end{equation}
in the vague topology on $(F,d_F)$. Finally, it is possible to check that $\nu(\partial B_F(x,r))=0$ for any $x\in F$ and $r>0$, exactly as in Lemma \ref{nucont}. Hence Assumption 1(c) does indeed hold.
\qed\end{proof}

We continue by considering the continuous time simple random walk $Y^\mathcal{G}$ on $\mathcal{G}$, which satisfies the following properties.

{\lem\label{carpetbounds} There exists a deterministic constant $d_w>2$ such that, $\mathbb{P}_\mu$-a.s., the transition density $\tilde{p}^\mathcal{G}$ of $Y^\mathcal{G}$ satisfies (\ref{nfghk}) and the associated exit time satisfies (\ref{exit}) for some constants $c_1,c_2,c_3,c_4\in(0,\infty)$. Furthermore, there $\mathbb{P}_\mu$-a.s. exist constants $c_5,c_6\in(0,\infty)$ such that
\begin{equation}\label{sclower}
\tilde{p}^\mathcal{G}_t(x,y)\geq c_5 t^{-d_f/d_w} \exp\left(-c_6\left(\frac{d_\mathcal{G}(x,y)^{d_w}}{t}\right)^{1/(d_w-1)}\right),\:\forall x,y\in V(\mathcal{G}),t\geq d_\mathcal{G}(x,y).
\end{equation}}
\begin{proof} In the unweighted case ($\mu_{xy}^\mathcal{G}=1$ for every $\{x,y\}\in E(\mathcal{G})$), the transition density bounds of (\ref{nfghk}) and (\ref{sclower}) are \cite{RWSC}, Theorem 7.1. The exit time bound of (\ref{exit}) can be proved in a similar way to \cite{RWSC}, Theorem 5.5. The results for the random case follow from these results using a rough isometry argument (for example, apply \cite{GT}, Theorem 3.1, and \cite{HamKum2}, Theorem 5.11).
\qed\end{proof}

The continuous time version of the tightness condition of Assumption 2 is an easy consequence of this result.

{\propn \label{carpetprop2} $\mathbb{P}_\mu$-a.s., the graphs $(G^n)_{n\geq 0}$ satisfy the continuous time version of Assumption 2 for $\alpha(n)$, $\beta(n)$ as in Proposition \ref{carpetprop} and $\gamma(n)=c_2L^{d_wn}$, for some deterministic constant $c_2\in(0,\infty)$.}
\begin{proof} The heat kernel bounds of (\ref{nfghk}) and (\ref{sclower}) imply that the continuous time version of the parabolic Harnack inequality with exponent $d_w$ holds for $\mathcal{G}$ (see \cite{GT}, Theorem 3.1, for the analogous discrete time result), and it follows that the continuous time version of Assumption 4 holds with $\kappa=d_w$. We can then apply the continuous time version of Proposition \ref{42} to deduce the proposition.
\qed\end{proof}

Unfortunately, for generalised Sierpinski carpets, a weak convergence result for $Y^\mathcal{G}$ has not yet been proved, even in the case of deterministic weights. However, the heat kernel bounds for the simple random walk do imply that if $\tilde{\mathbf{P}}^n_\rho$ is the law of $(L^{-n}Y^\mathcal{G}_{\gamma(n)t})_{t\geq 0}$ under $\tilde{\mathbf{P}}^\mathcal{G}_\rho$, considered as a probability measure on $D([0,1], F)$, then the sequence $(\tilde{\mathbf{P}}^n_\rho)_{n\geq 0}$ is tight. Consequently, it holds that $(\tilde{\mathbf{P}}^n_\rho)_{n\geq 0}$ admits a convergent subsequence, and we will later show that for any convergent subsequence there is a corresponding local limit theorem.

{\lem $\mathbb{P}_\mu$-a.s., the sequence $(\tilde{\mathbf{P}}^n_\rho)_{n\geq 0}$ is tight in $D(\mathbb{R}_+, F)$ and, moreover, if $\mathbf{P}_\rho$ is a limit point of $(\tilde{\mathbf{P}}^n_\rho)_{n\geq 0}$, then $\mathbf{P}_\rho(C(\mathbb{R}_+,F))=1$.}
\begin{proof} Let $t\geq 0$ and $\varepsilon>0$, then (\ref{exit}) implies that, $\mathbb{P}_\mu$-a.s.,
\[\lim_{\delta\rightarrow 0}\limsup_{n\rightarrow\infty} \delta^{-1}\tilde{\mathbf{P}}^\mathcal{G}_\rho \left(\sup_{s\in[t,t+\delta]}\left|L^{-n}Y^\mathcal{G}_{\gamma(n)s}-L^{-n}Y^\mathcal{G}_{\gamma(n)t}\right|\geq \varepsilon\right)=0.\]
The lemma follows (cf. Theorem 7.3 and the corollary to Theorem 7.4 in \cite{Bill2}).
\qed\end{proof}

In view of this result and Propositions \ref{carpetprop} and \ref{carpetprop2}, to apply Theorem \ref{bdd} to deduce local limit theorems along convergent subsequences of $(\tilde{\mathbf{P}}^n_\rho)_{n\geq 0}$, it suffices to show that any limit point, $\mathbf{P}_\rho$ say, admits a family of transition densities $(q_t(x))_{x\in F,t>0}$ that is jointly continuous in $(t,x)$. To prove that this is the case, we first extend Proposition \ref{carpetprop2}. We write $\tilde{p}^n_t(x,y)=\tilde{p}^{{G}^n}_t(x,y)$ and $\tilde{q}^n_t(x)=\tilde{p}^{{G}^n}_t(\rho,x)$.

{\lem \label{tightcarpet} $\mathbb{P}_\mu$-a.s., for any compact interval $I\subset (0,\infty)$ and $r>0$,
\[\lim_{\delta\rightarrow 0}\limsup_{n\rightarrow\infty}\sup_{\substack{x,y\in B_{G^n}(\rho,\alpha(n)r):\\d_{G^n}(x,y)\leq \alpha(n)\delta}}\sup_{\substack{s,t \in I:\\|s-t|\leq\delta}}\beta(n)\left|\tilde{q}_{\gamma(n)s}^{n}(x)-\tilde{q}_{\gamma(n)t}^{n}(y)\right|=0.\]}
\begin{proof} Given Proposition \ref{carpetprop2}, it is enough to demonstrate that $\mathbb{P}_\mu$-a.s. for any compact interval $I\subset (0,\infty)$ and $r>0$,
\begin{equation}\label{aim}
\lim_{\delta\rightarrow 0}\limsup_{n\rightarrow\infty}\sup_{x\in B_{G^n}(\rho,\alpha(n)r)}\sup_{\substack{s,t \in I:\\|s-t|\leq\delta}}\beta(n)\left|\tilde{q}_{\gamma(n)s}^{n}(x)-\tilde{q}_{\gamma(n)t}^{n}(x)\right|=0.
\end{equation}
The following argument holds $\mathbb{P}_\mu$-a.s. We can write, for $s<t$,
\begin{eqnarray*}
\beta(n)\left|\tilde{q}_{\gamma(n)s}^{n}(x)-\tilde{q}_{\gamma(n)t}^{n}(x)\right|
&\leq &\beta(n)\int_{F}\left|\tilde{q}_{\gamma(n)s}^{n}(x)-\tilde{q}_{\gamma(n)s}^{n}(y)\right|\tilde{p}^n_{\gamma(n)(t-s)}(y,x)\nu^n(dy),
\end{eqnarray*}
where $\nu^n:=\nu^{G^n}$. Now, for $\eta>0$,
\begin{eqnarray*}
\lefteqn{\lim_{\delta\rightarrow 0}\limsup_{n\rightarrow\infty}\sup_{x\in B_{G^n}(\rho,\alpha(n)r)}\sup_{\substack{s,t \in I:\\|s-t|\leq\delta}}\beta(n)\int_{B_F(x,\eta)}\left|\tilde{q}_{\gamma(n)s}^{n}(x)-\tilde{q}_{\gamma(n)s}^{n}(y)\right|}\\
&&\hspace{200pt}\times\tilde{p}^n_{\gamma(n)(t-s)}(y,x)\nu^n(dy)\\
&\leq&\limsup_{n\rightarrow\infty}\sup_{x\in B_{G^n}(\rho,\alpha(n)r)}\sup_{y\in B_F(x,\eta)}\sup_{s\in I} \beta(n)\left|\tilde{q}_{\gamma(n)s}^{n}(x)-\tilde{q}_{\gamma(n)s}^{n}(y)\right|.\hspace{110pt}
\end{eqnarray*}
For any given $\varepsilon>0$, by Proposition \ref{carpetprop2}, we can make this upper bound smaller than $\varepsilon$ by choosing $\eta$ small enough. Fixing such an $\eta$, we also have
\begin{eqnarray*}
\lefteqn{\lim_{\delta\rightarrow 0}\limsup_{n\rightarrow\infty}\sup_{x\in B_{G^n}(\rho,\alpha(n)r)}\sup_{\substack{s,t \in I:\\|s-t|\leq\delta}}\beta(n)\int_{F\backslash B_F(x,\eta)}\left|\tilde{q}_{\gamma(n)s}^{n}(x)-\tilde{q}_{\gamma(n)s}^{n}(y)\right|}\\
&&\hspace{200pt}\times\tilde{p}^n_{\gamma(n)(t-s)}(y,x)\nu^n(dy),\\
&\leq&\lim_{\delta\rightarrow 0}\limsup_{n\rightarrow\infty}\left(\sup_{x\in V(\mathcal{G})}\sup_{s\in I} 2\beta(n)\tilde{q}_{\gamma(n)s}^{n}(x)\right)\left(\sup_{x\in V(\mathcal{G})}
\tilde{\mathbf{P}}^\mathcal{G}_x\left(\tau(x,\alpha(n)\eta)\leq \gamma(n)\delta\right)\right)\hspace{30pt}\\
&=&0,
\end{eqnarray*}
where we apply the bounds of Lemma \ref{carpetbounds} to deduce the final equality. The limit result at (\ref{aim}) follows.
\qed\end{proof}

{\lem $\mathbb{P}_\mu$-a.s., if $\mathbf{P}_\rho$ is a limit point of the sequence $(\tilde{\mathbf{P}}^n_\rho)_{n\geq 0}$, then $\mathbf{P}_\rho$ admits a family of transition densities $(q_t(x))_{x\in F,t>0}$ that is jointly continuous in $(t,x)$.}
\begin{proof} In this proof, which holds $\mathbb{P}_\mu$-a.s., we fix a subsequence $(n_i)_{i\geq 0}$ so that $(\tilde{\mathbf{P}}^{n_i}_\rho)_{i\geq 0}$ converges, and let $\mathbf{P}_\rho$ be the corresponding limit point. Applying Lemmas \ref{carpetmetric} and \ref{tightcarpet}, it is elementary to define, for every $n\geq 0$, a jointly continuous function $(f^n(t,x))_{x\in F,t>0}$, such that $f^n(t,x)=\beta(n)\tilde{q}^n_{\gamma(n)t}(x)$ for every $x\in V(G^n)$, $t>0$, in such a way that the sequence $(f^n)_{n\geq 0}$ is tight in $C(I\times \overline{B}_F(\rho,r),\mathbb{R})$ for any compact interval $I\subset (0,\infty)$ and $r>0$.

Fix a compact interval $I\subset (0,\infty)$ and choose $r$ large enough so that $K=E_0\cap F$ is contained inside ${B}_F(\rho,r)$. By the conclusion of the previous paragraph, it is possible to choose a subsequence $(n_{i_j})_{j\geq 0}$ and jointly continuous $f^{I,r} =(f^{I,r}(t,x))_{x\in \overline{B}_F(\rho,r) ,t\in I}$ such that $f^{I,r}$ is the uniform limit of $(f^{n_{i_j}})_{j\geq 0}$ on $I\times \overline{B}_F(\rho,r)$. Now, suppose $S\in \mathcal{S}_n(E_0)$, for some $n\geq 0$, and write $A=S\cap F$. If we define $A_\varepsilon:=\{x\in F:\:d_F(x,A)<\varepsilon\}$, then $A_\varepsilon$ is an open subset of $F$, and consequently (see \cite{Bill2}, Theorem 2.1, for example), for $t>0$,
\[\liminf_{j\rightarrow\infty}\tilde{\mathbf{P}}_\rho^{n_{i_j}}(X_t\in A_\varepsilon)\geq \mathbf{P}_\rho(X_t\in A_\varepsilon).\]
Moreover, by the definitions of $(f^n)_{n\geq 0}$ and $f^{I,r}$, if $t\in I$,
\[\liminf_{j\rightarrow\infty}\tilde{\mathbf{P}}_\rho^{n_{i_j}}(X_t\in A_\varepsilon)\leq\limsup_{j\rightarrow\infty}\int_{\overline{A}_{\varepsilon}}f^{n_{i_j}}(t,x)\beta(n_{i_j})^{-1}\nu^{n_{i_j}}(dx)\leq \int_{\overline{A}_{\varepsilon}}f^{I,r}(t,x)\nu(dx),\]
where $\overline{A}_{\varepsilon}$ is the closure of $A_\varepsilon$ and we also apply the measure convergence of (\ref{vaguecarpet}). Letting $\varepsilon\rightarrow 0$, we obtain $\int_{{A}}f^{I,r}(t,x)\nu(dx)\geq \mathbf{P}_\rho(X_t\in A)$. Similarly, by first considering $A_{-\varepsilon}:=\{x\in A:\:d_F(x,\partial A)\geq \varepsilon\}$, we are able to prove that the opposite inequality is also true. An elementary $\sigma$-algebra argument (\cite{Kallenberg}, Lemma 1.17, for example) allows this result to be extended to show that
\begin{equation}\label{td}
\int_{{A}}f^{I,r}(t,x)\nu(dx)= \mathbf{P}_\rho(X_t\in A),
\end{equation}
for every measurable $A\subseteq K$, $t\in I$. Noting that $\nu(A)>0$ for every open set $A$, it follows that the choice of subsequence is unimportant and $f^{I,r}$ is actually the uniform limit of $(f^{n_i})_{i\geq 0}$ in the region $I\times K$. Repeating the same argument on an increasing sequence of space-time regions, we can extend the definition of $f^{I,r}$ to deduce the existence of a jointly continuous function $f=f(t,x)_{x\in F,t>0}$ that is the point-wise limit of $(f^{n_i})_{i\geq 0}$ everywhere in $(0,\infty)\times F$, with the limit being uniform on compacts, and, moreover, (\ref{td}) holds with $f^{I,r}$ replaced by $f$ for any measurable $A\subseteq F$ and $t>0$.
\qed\end{proof}

We now can state the main conclusion of this section, which is an application of Theorem \ref{bdd}. Note that the heat kernel bounds of (\ref{nfghk}) and (\ref{sclower}) imply that the densities defined in the previous lemma satisfy $q_t(x)\neq 0$ for every $x\in F$ and $t>0$, so the limit is non-trivial.

{\thm For $\mathbb{P}_\mu$-a.e. realisation of a Sierpinski carpet graph $\mathcal{G}$ with independent and identically-distributed edge-weights that satisfy uniform boundedness: if we fix a compact interval $I\subset (0,\infty)$ and $r>0$, suppose that $(\tilde{\mathbf{P}}^{n_i}_\rho)_{i\geq 0}$ converges, and let $\mathbf{P}_\rho$ and $(q_t(x))_{x\in F,t>0}$ represent the corresponding limit point and family of transition densities, then
\[\lim_{i\rightarrow\infty} \sup_{x,y\in B_F(\rho,R)}\sup_{t\in I}\left|c_1N^{n_i}\tilde{q}^\mathcal{G}_{c_2L^{d_wn_i} t}(\tilde{g}_{n_i}(x))-q_t(x)\right|=0.\]}

Finally, note that if (a subsequence of) $(\tilde{\mathbf{P}}^n_\rho)_{n\geq 0}$ was shown to converge to the Brownian motion on the Sierpinski carpet, as constructed in \cite{BarBashd} (see \cite{BarBas} for the two-dimensional case), then the transition density estimates of Lemma \ref{carpetbounds} and \cite{BarBashd}, Theorem 1.3, (cf. \cite{BarBastd}, Theorem 1.1), would enable us to apply Theorem \ref{unbdd} to extend the above result to unbounded regions of time and space.

\def\cprime{$'$}
\providecommand{\bysame}{\leavevmode\hbox to3em{\hrulefill}\thinspace}
\providecommand{\MR}{\relax\ifhmode\unskip\space\fi MR }
\providecommand{\MRhref}[2]{%
  \href{http://www.ams.org/mathscinet-getitem?mr=#1}{#2}
}
\providecommand{\href}[2]{#2}

\end{document}